\documentstyle[12pt]{article}
 \input amssym.def
\input amssym
\begin{document}
\def \Z{\Bbb Z}
\def \C{\Bbb C}
\def \R{\Bbb R}
\def \Q{\Bbb Q}
\def \N{\Bbb N}
\def \wt{{\rm wt}}
\def \tr{{\rm tr}}
\def \span{{\rm span}}
\def \Res{{\rm Res}}
\def \Res{{\rm QRes}}
\def \End{{\rm End}}
\def \E{{\rm End}}
\def \Ind {{\rm Ind}}
\def \Irr {{\rm Irr}}
\def \Aut{{\rm Aut}}
\def \Hom{{\rm Hom}}
\def \mod{{\rm mod}}
\def \ann{{\rm Ann}}
\def \<{\langle} 
\def \>{\rangle} 
\def \t{\tau }
\def \a{\alpha }
\def \e{\epsilon }
\def \l{\lambda }
\def \L{\Lambda }
\def \g{\gamma}
\def \b{\beta }
\def \om{\omega }
\def \o{\omega }
\def \c{\chi}
\def \ch{\chi}
\def \cg{\chi_g}
\def \ag{\alpha_g}
\def \ah{\alpha_h}
\def \ph{\psi_h}
\def \be{\begin{equation}\label}
\def \ee{\end{equation}}
\def \bl{\begin{lem}\label}
\def \el{\end{lem}}
\def \bt{\begin{thm}\label}
\def \et{\end{thm}}
\def \bp{\begin{prop}\label}
\def \ep{\end{prop}}
\def \br{\begin{rem}\label}
\def \er{\end{rem}}
\def \bc{\begin{coro}\label}
\def \ec{\end{coro}}
\def \bd{\begin{de}\label}
\def \ed{\end{de}}
\def \pf{{\bf Proof. }}
\def \voa{{vertex operator algebra}}

\newtheorem{thm}{Theorem}[section]
\newtheorem{prop}[thm]{Proposition}
\newtheorem{coro}[thm]{Corollary}
\newtheorem{conj}[thm]{Conjecture}
\newtheorem{lem}[thm]{Lemma}
\newtheorem{rem}[thm]{Remark}
\newtheorem{de}[thm]{Definition}
\newtheorem{hy}[thm]{Hypothesis}
\makeatletter
\@addtoreset{equation}{section}
\def\theequation{\thesection.\arabic{equation}}
\makeatother
\makeatletter

\baselineskip=16pt
\begin{center}{\Large \bf Hom functor
 and the associativity of tensor products of modules for vertex operator
algebras}

\vspace{0.5cm}
Chongying Dong\footnote{Supported by NSF grant DMS-9303374 and a
research grant from the Committee on Research, UC Santa Cruz.},
Haisheng Li and Geoffrey Mason\footnote{Supported by NSF grant
DMS-9401272 and a research grant from the Committee on Research, UC
Santa Cruz.}\\ Department of Mathematics, University of California,
Santa Cruz, CA 95064
\end{center}

\begin{abstract} We give a new, construction-free proof of the
associativity of tensor product for modules for rational vertex
operator algebras under certain convergence conditions.
\end{abstract}

\section{Introduction}

Vertex operator algebras (or chiral algebras as
commonly used in physical references) are closely related to classical
mathematical notions such as Lie algebras and associative algebras.
The study of tensor products of representations for a vertex operator
algebra was initiated by physicists in order to explain  
fusion rules (cf. [MS]). Mathematically, a
tensor product theory for modules for a vertex operator algebra has
been developed in [HL0]-[HL4] and [KL0]-[KL2], so that the classical
Clebsch-Gordon coefficients give the fusion rules.  

In classical Lie theory, both the existence and associativity of the
tensor product of two modules follow from the Hopf algebra structure
of the universal enveloping algebra of a Lie algebra. However, the
tensor product theory for vertex operator algebras is much more
complicated.  Especially, the proof of the associativity is highly
nontrivial.  

Using the tensor product construction given in [HL0]-[HL4], 
Huang has established in [H1] the associativity of tensor
products under a certain assumption on convergence of products of two
intertwining operators. In the present paper we will give an alternate 
approach to Huang's result which is based on the definition of tensor 
product in terms of certain universal properties as developed in [LL]
and [L1]. Our formal variable approach yields a proof that is 
shorter and simpler than that of Huang, and has the advantage of being
manifestly {\em construction-free} in the sense that it does not
depend on the details of the construction of the tensor product.
(Huang has pointed out to us that the proof in [H1] may be adapted
so that it too is construction-free in this sense.) 

Here we explain our approach by using the classical tensor product
theory for a Lie algebra $L$. The reader will find that our proof of
associativity of the tensor product for modules for a vertex operator algebra
is similar. 
For three $L$-modules $U_{i}$
$(i=1,2,3)$, an intertwining operator of type
$\left(\!\begin{array}{c}U_{3}\\U_{1}U_{2}\end{array}\!\right)$ is a
linear map $f$ from $U_{1}$ to ${\rm Hom}_{{\C}}(U_{2},U_{3})$ such
that $af(u_{1})u_{2}=f(au_{1})u_2+f(u_{1})au_{2}$ for $a\in L, u_{i}\in
U_{i}$.  A tensor product for $(U_{1},U_{2})$ can be defined as a pair
$(U, f)$, where $U$ is a $L$-module and $f$ is an intertwining
operator of type
$\left(\!\begin{array}{c}U\\U_{1}U_{2}\end{array}\!\right)$ satisfying
the following universal property: for any $L$-module $W$ and for any
intertwining operator $g$ of type
$\left(\!\begin{array}{c}W\\U_{1}U_{2}\end{array}\!\right)$, there
exists a unique $L$-homomorphism $\psi$ from $U$ to $W$ such that
$g=\psi f$.

To avoid an explicit construction of the tensor product module, we may
interpret $U_{1}\otimes (U_{2}\otimes U_{3})$ as a quadruple
$(W_{1,23},W_{2,3}, f_{1,23},f_{2,3})$, where $W_{1,23}$ and $W_{2,3}$
are $L$-modules and $f_{1,23}$ and $f_{2,3}$ are intertwining operator
of types
$\left(\begin{array}{c}W_{1,23}\\U_{1}W_{2,3}\end{array}\right)$ and
$\left(\begin{array}{c}W_{2,3}\\U_{2}U_{3}\end{array}\right)$
respectively. The quadruple satisfies the following universal property: for
any quadruple $(W_{1},W_{2}, f_{1},f_{2})$, where $W_{i}$ are
$L$-modules and $f_{1}$ is an intertwining operator of type
$\left(\begin{array}{c}W_{1}\\U_{1}W_{2} \end{array}\right)$ and
$f_{2}$ is an intertwining operator of type
$\left(\begin{array}{c}W_{2}\\U_{2}U_{3}\end{array}\right)$, 
there exists a unique $L$-homomorphism $\psi$ from $W_{1,23}$ to
$W_{1}$ such that $$f_{1}(u_{1})f_{2}(u_{2})u_{3}=\psi
f_{1,23}(u_{1})f_{2,3}(u_{2})u_{3}$$ for $u_{i}\in U_{i}$.  Similarly,
we  may also interpret the tensor product module $(U_{1}\otimes
U_{2})\otimes U_{3}$ as a quadruple $(W_{12,3},W_{1,2},
f_{12,3},f_{1,2})$.

In order to see that $W_{1,23}$ and $W_{12,3}$ are isomorphic
$L$-modules, we construct 
two $L$-homomorphisms (between $W_{1,23}$ and
 $W_{12,3}$) which are inverses to each other 
by using the universal properties for both $W_{1,23}$ and $W_{12,3}$.  
For any $u_{i}\in U_{i}$ $(i=2,3)$, we obtain an element $F(u_{2},u_{3})\in
{\rm Hom}_{{\C}}(U_{1},W_{12,3})$ such that
$F(u_{2},u_{3})u_{1}=f_{12,3}(f_{1,2}(u_{1})u_{2})(u_{3})$ for any
$u_{1}\in U_{1}$.  Then all these elements form a submodule of ${\rm
Hom}_{{\C}}(U_{1},W_{12,3})$.  Denote this module by $W_{2}$. Then we
have an intertwining operator $f_{2}$ of type
$\left(\!\begin{array}{c}W_{2}\\U_{2}U_{3}\end{array}\!\right)$
defined by $f_{2}(u_{2})u_{3}=F(u_{2},u_{3})$ for $u_{i}\in
U_{i}$. Furthermore, we obtain an intertwining operator $f_{1}$ of
type
$\left(\!\begin{array}{c}W_{12,3}\\U_{1}W_{2}\end{array}\!\right)$
defined by $f_{1}(f_{2}(u_{2})u_{3})u_{1}=F(u_{2},u_{3})u_{1}$ for
$u_{i}\in U_{i}$.  Then it follows from the universal property for
$W_{1,23}$ that there exists a unique $L$-homomorphism $\psi_{1}$ from
$W_{1,23}$ to $W_{12,3}$ such that
$$f_{1}(u_{1})f_{2}(u_{2})u_{3}=\psi
f_{1,23}(u_{1})f_{2,3}(u_{2})u_{3}$$ for $u_{i}\in U_{i}$.  Similarly,
we obtain an $L$-homomorphism $\psi_{2}$ from $W_{12,3}$ to
$W_{1,23}$.  Then it easily follows that $W_{1,23}$ and $W_{12,3}$ are
isomorphic $L$-modules.

This paper consists of three sections. We review the analogue of
``Hom''-functor and its relation to tensor products of modules for a vertex
operator algebra from [L1] and [LL]  in Section 2. Section 3 is devoted
to the proof of associativity for modules.

\section{Tensor product and ``Hom''-functor}
In this section we recall some basic  definitions such as (weak) 
module, rationality, intertwining 
operator and tensor product and we
reformulate some results in [LL]. 
{\em Throughout the paper, $V$ will be a vertex operator algebra.}
The reader is referred to [B], [FLM], [FHL] for the definition of vertex
operator algebra.

We use $z, z_{0}, z_{1}, \cdots$ for commuting formal variables.

A {\em weak $V$-module} is a 
pair $(W,Y_{W})$, where $W$ is a vector space and $Y_{W}(\cdot,z)$ is a linear map from 
$V$ to $({\rm End} W)[[z,z^{-1}]]$ 
satisfying the following axioms:
(1) $Y_{W}({\bf 1},z)=id_{W}$; (2) $Y_{W}(a,z)u\in W((z))$ for any $a\in V,u\in W$;
(3) $Y_{W}(L(-1)a,z)={d\over dz}Y_{W}(a,z)$ for $a\in V$; (4) the Jacobi identity:
\begin{eqnarray}
& &z_{0}^{-1}\delta\left(\frac{z_{1}-z_{2}}{z_{0}}\right)Y_{M}(a,z_{1})Y_{M}(b,z_{2})u
-z_{0}^{-1}\delta\left(\frac{z_{2}-z_{1}}{-z_{0}}\right)Y_{M}(b,z_{2})Y_{M}(a,z_{1})u
\nonumber\\
& &=z_{2}^{-1}\delta\left(\frac{z_{1}-z_{0}}{z_{2}}\right)Y_{M}(Y(a,z_{0})b,z_{2})u
\end{eqnarray}
for $a,b\in V, u\in W$. As observed in [DLM1], axiom (3) is redundant.
We define ${\cal{C}}^{0}(V)$ to be the category of all weak $V$-modules
(with the obvious morphisms).

A weak $V$-module $(W,Y_{W})$ 
is called a {\em $V$-module} if $L(0)$ acts semisimply on $W$ with 
 $L(0)$-eigenspace decomposition
$M=\oplus_{h\in {\C}}M_{h}$ such that for any $h\in {\C}$,
$\dim M_{h}<\infty$ and $M_{h+n}=0$ for  $n\in {\Z}$ sufficiently small.
We denote by ${\cal{C}}(V)$ the category of all $V$-modules. It is
a full subcategory of ${\cal C}^0(V).$  

For a weak $V$-module $W,$ $0\ne w\in W$ is called a {\em generalized
weight vector with generalized weight} $h$ if there exists some $n$
such that $(L(0)-h)^nw=0.$ We write gwt$w=h$ in this case and call $w$
a {\em homogeneous} vector. We denote by $P(W)$ the set of all
generalized weights of $W$.  Let ${\cal{C}}^{1}(V)$ be the category of
weak $V$-modules $W$ such that $L(0)$ acts locally finitely on $W$
(that is, for every $w\in W$ there exists a finite-dimensional
subspace $M$ of $W$ which contains $w$ and is $L(0)$-invariant), and
such that $P(W)\subseteq \cup_{i=1}^{n}(h_{i}+\Z_{+})$ for finitely
many complex numbers $h_{1},\cdots, h_{n}$. Note that in this case $W$
is a direct sum of generalized eigenspaces for $L(0).$

A {\em ${\Z}_{+}$-graded weak} $V$-module ([FZ], [DLM1])
 is  a weak $V$-module $W$
together with a ${\Z}_{+}$-grading
$W=\oplus_{n=0}^{\infty}W(n)$  such that
\begin{eqnarray}
a_{m}W(n)\subseteq W(k+n-m-1)\;\;\;\mbox{ for }a\in V_{(k)}, m,n\in {\Z},
\end{eqnarray}
where $W(n)=0$ by convention if $n<0$. Denote by ${\cal{C}}^{+}(V)$ the
category of $\Z_{+}$-graded weak $V$-modules. One can prove that 
${\cal{C}}^{1}(V)$
is a subcategory of ${\cal{C}}^{+}(V)$.

One may define the notions of ``submodule''and 
``irreducible submodule'' accordingly. A VOA $V$ is said to be {\em rational}
if the category ${\cal{C}}^{+}(V)$ is semisimple, {\it i.e.}, 
any ${\Z}_{+}$-graded weak $V$-module is a direct sum of irreducible ${\Z}_{+}$-graded 
weak $V$-modules. It was proved in [DLM1] that if $V$ is rational then
there are only finitely many
irreducible ${\Z}_{+}$-graded weak $V$-modules up to equivalence and any irreducible weak 
$V$-module is a module. Consequently, this definition of rationality agrees with Zhu's 
definition of rationality [Z]. If $V$ is rational, then 
${\cal{C}}^{+}(V)={\cal{C}}^{1}(V)$ contains ${\cal{C}}(V)$ as a subcategory and 
they have the same irreducible objects.
It was proved in [DLM2] that
${\cal{C}}^{0}(V)={\cal{C}}^{+}(V)={\cal{C}}^{1}(V)$ for almost all known rational 
vertex operator algebras such as $V_{L}$, associated with a positive definite even 
lattice $L$; $L(\ell,0)$,
associated with an affine Lie algebra with a positive integral level $\ell$; 
$L(c_{p,q},0)$, associated with the Virasoro algebra with a central charge
$c_{p,q}$ in the minimal series; and $V^{\natural},$ the moonshine module.

Let $W_{i}$ $(i=1,2,3)$ be weak $V$-modules. {\em An intertwining operator} of type
$\left(\begin{array}{c}W_{3}\\W_{1}W_{2}\end{array}\right),$
as defined in [FHL], is a linear map $I(\cdot,z)$ from $W_{1}$ to 
$({\rm Hom}_{\C}(W_{2},W_{3}))\{z\}$
satisfying  the truncation condition (2), the $L(-1)$-derivative property (3), and 
the Jacobi 
identity (4), where  $Y_{M}$ is replaced by $I$.

\bd{d2.1} 
Let ${\cal{D}}$ be any one of the categories 
${\cal{C}}^{0}(V), {\cal{C}}^{1}(V), {\cal{C}}^{+}(V)$ or ${\cal{C}}(V)$.
Let $W_{1}$ and $W_{2}$ be two weak $V$-modules from the category ${\cal{D}}$. A
{\em tensor product} for the ordered pair $(W_{1},W_{2})$ is a pair
$(M,F(\cdot,z))$ consisting of a weak $V$-module $M$ in the category $\cal{D}$
and an intertwining
operator $F(\cdot,z)$ of type
$\left(\!\begin{array}{c}M\\W_{1} W_{2}\end{array}\!\right)$
satisfying the following universal 
property: For any weak $V$-module $W$ in the category $\cal{D}$ and 
any intertwining
operator $I(\cdot,z)$ of type
$\left(\!\begin{array}{c}W\\W_{1} W_{2}\end{array}\!\right)$, there
exists a unique $V$-homomorphism $\psi$ from $M$ to $W$ such that
$I(\cdot,z)=\psi\circ F(\cdot,z)$. (Here $\psi$ extends canonically to a
linear map from $M\{z\}$ to $W\{z\}$.) We shall denote the tensor product
$M$ by $W_1\boxtimes_{\cal D}W_2.$ 
\ed

\begin{lem}\label{l2.3}
Let $(W,F(\cdot,z))$ be a tensor product for the
ordered pair $(W_{1},W_{2})$. Then $F(\cdot,z)$ is surjective in the
sense that all the coefficients of the formal series $F(u_{1},z)u_{2}$ for $u_{i}\in
W_{i}$ linearly span $W$.
\end{lem}

\begin{coro}\label{c2.4} 
If $(M,F(\cdot,z))$ is a tensor product for the
ordered pair $(W_{1},W_{2})$ of weak $V$-modules in the category ${\cal{D}}$, then 
for any weak $V$-module
$W_{3}$ in $\cal{D}$, ${\rm Hom}_{V}(M,W_{3})$ is linearly isomorphic 
to the space
of intertwining operators of type
$\left(\!\begin{array}{c}W_{3}\\W_{1} W_{2}\end{array}\!\right)$.
\end{coro}

\begin{rem}\label{r2.2} 
Just as in the classical algebra theory, it follows
from the universal property that if there exists a tensor product
$(M,F(\cdot,z))$
for the ordered pair $(W_{1},W_{2})$, then it is unique up to 
$V$-isomorphism, i.e., if $(W,G(\cdot,z))$ is another tensor
product, then there is a $V$-isomorphism $\psi$ from $M$ to $W$
such that $G=\psi \circ F$. Conversely, let $(M,F(\cdot,z))$ be a tensor
product for the ordered pair $(W_{1},W_{2})$ and let $\sigma$ be a 
$V$-isomorphism from $M$ to $W$. Then $(W, \sigma\circ F(\cdot,z))$
is a tensor product for $(W_{1},W_{2})$.
\end{rem}

\begin{rem}\label{r2.3}
Assuming that $V$ is a rational vertex operator algebra such that the
 fusion rules among any three irreducible $V$-modules are finite, the
 existence of a tensor product in the category ${\cal{C}}(V)$ for any
 pair $(W_{1},W_{2})$ has been proved in [HL1]-[HL4] and [LL]. On the
 other hand, it follows from [LL] that tensor products exist in the
 category ${\cal{C}}^{+}(V)$ without assuming the finiteness of fusion
 rules.
\end{rem}

Sometimes, even if $V$ is not rational, tensor products for certain
categories may still exist.  For instance, let $\frak{g}$ be a simple
Lie algebra with dual Coxeter number $\Omega$ and let $L(\ell,0)$ be the
vertex operator algebra associated with the affine Lie algebra
$\hat{\frak{g}}$ of level $\ell$.  If $\ell+\Omega \notin \Q_{+}$,
then $L(\ell,0)$ is not rational.  In [KL0]-[KL2], tensor products for
a certain subcategory of the category ${\cal{C}}^{1}(V)$ have been
constructed, although the whole setting is at the level of affine Lie algebras
rather than vertex operator algebras.

\begin{rem}\label{r2.4}
Let $V$ be a vertex operator algebra. Then 
any weak $V$-module $W$ from ${\cal{C}}^{1}(V)$ has a canonical submodule decomposition
$W=W_{1}\oplus \cdots \oplus W_{n}$ such that $P(W_{i})\subseteq h_{i}+\Z_{+}$ 
for $i=1,\cdots, n$.
For any $\lambda\in \C$, we define ${\cal{C}}_{\lambda}(V)$ to be the subcategory
of ${\cal{C}}^{1}(V)$ with objects $W$ satisfying $P(W)\subseteq \lambda+\Z_{+}$. 

Let $W_{1}$ and $W_{2}$ be objects in the categories
${\cal{C}}_{\lambda_1}(V)$ and ${\cal{C}}_{\lambda_2}(V)$,
respectively.  {}From Theorem 5.2.17 [L2], for any $\lambda\in \C$ we
have a weak $V$-module $T(W_{1},W_{2})^{\lambda}$ in the category
${\cal{C}}_{\lambda}(V)$ and an intertwining operator
$F_{\lambda}(\cdot,z)$ of type
$\left(\!\begin{array}{c}T(W_{1},W_{2})^{\lambda}\\W_{1}W_{2}\end{array}\!\right)$
satisfying the following universal property: for any $W_{3}$ from
${\cal{C}}_{\lambda}(V)$ and any intertwining operator $I(\cdot,z)$ of
type $\left(\!\begin{array}{c}W_{3}\\W_{1}W_{2}\end{array}\!\right)$,
there exists a unique $V$-homomorphism $f$ from
$T(W_{1},W_{2})^{\lambda}$ to $W_{3}$ such that
$I(\cdot,z)=fF_{\lambda}(\cdot,z)$. One may consider
$T(W_{1},W_{2})^{\lambda}$ as a product of $W_{1}$ and $W_{2}$ in the
category ${\cal{C}}_{\lambda}(V)$.  Furthermore, there is a
$V$-homomorphism $f_{\lambda}$ from $T(W_{1},W_{2})^{\lambda}$ onto
$T(W_{1},W_{2})^{\lambda-1}$ such that
$F_{\lambda-1}(\cdot,z)=f_{\lambda}F_{\lambda}(\cdot,z)$.

Set
$$P(W_{1},W_{2})=\cup_{W\in Irr(V), f(W_{1},W_{2};W)\ne 0}P(W),$$
where $f(W_{1},W_{2};W)$ is the fusion rule of type
$\left(\!\begin{array}{c}W\\W_{1}W_{2}\end{array}\!\right)$.

Assume that there are finitely many complex numbers $h_{1},\cdots, h_{k}$ such that
$P(W_{1},W_{2})$ $\subseteq \cup_{i=1}^{k}(h_{i}+\Z_{+}),$ and that 
$h_{1},\cdots, h_{k}$ are distinct modulo $\Z$.
Then one can
show that the direct sum of the $T(W_{1},W_{2})^{h_{i}}$ for $i=1,\cdots,k$ is
a tensor product of $W_{1}$ with $W_{2}$ in the category 
${\cal{C}}^{1}(V)$. Since ${\cal{C}}^{1}(V)$ is closed under the operation of finite sum, 
it is clear that
tensor products exist for any pair of weak modules from ${\cal{C}}^{1}(V)$ under 
the above assumption.
\end{rem}

For our purposes we need to reformulate the notion of generalized
intertwining operator and some related results as follows:

\bd{d3.1} Let $W_{1}$ and $W_{2}$ be weak 
$V$-modules in a category ${\cal{D}}$ of (weak) $V$-modules.  An element 
$\Phi (z)=\displaystyle{\sum_{\alpha\in {\C}}\Phi_{\alpha}z^{-\alpha-1}}
\in \left({\rm Hom}_{{\C}}(W_{1},W_{2})\right)\{z\}$ is called a
{\it generalized left intertwining operator} from $W_{1}$ 
to $W_{2}$ if the following
conditions (Gl1)-(Gl3) hold: 

(Gl1)\hspace{0.25cm}  For any $\alpha\in {\C}, u_{1}\in W_{1}$,
$\Phi_{\alpha+n}u_{1}=0$ for $n\in {\Z}$ sufficiently large;

(Gl2)\hspace{0.25cm} $[L(-1),\Phi (z)]=\Phi'(z)
 \left(=\displaystyle{{d\over dz}}\Phi(z)\right);$

(Gl3)\hspace{0.25cm} For any $a\in V$, there exists a positive integer
$k$ such that 
\begin{eqnarray}\label{el}
(z_{1}-z_{2})^{k}Y_{W_{2}}(a,z_{1})\Phi(z_{2})
=(z_{1}-z_{2})^{k}\Phi(z_{2})Y_{W_{1}}(a,z_{1}).
\end{eqnarray}
\ed

Denote by $G_{l}(W_{1},W_{2})$ the space of all generalized left
intertwining operators from $W_{1}$ to $W_{2}$.  Suppose that $W_{1}$
and $W_{2}$ are objects in the category ${\cal{C}}^{1}(V)$.  A
homomorphism $f$ from $W_{1}$ to $W_{2}$ is said to be {\it
homogeneous of generalized weight $h$} if $f((W_{1})_{h_{1}})\subseteq
(W_{2})_{h+h_{1}}$ for $h_{1}\in \C$.  A generalized left intertwining
operator $\Phi(z)=\sum_{n\in \C}\Phi_{n}z^{-n-1}$ is said to be {\it
homogeneous of generalized weight $h$} if for any $n\in \C$,
$\Phi_{n}$ is a homogeneous operator of generalized weight $h-n-1$.  A
sufficient condition for this to occur is that $\Phi(z)$ satisfies
the following equation:
\begin{eqnarray}\label{elh}
[L(0),\Phi(z)]=\left(h+z{d\over dz}\right)\Phi(z). 
\end{eqnarray}
A generalized left intertwining operator $\Phi(z)$ of weight $h$ is said to 
be {\em primary} if the following condition holds:
\begin{eqnarray}\label{ep2}
[L(m), \Phi(z)]=z^{m}\left((m+1)h+z{d\over dz}\right)\Phi(z)\;\;\;\mbox{ for }m\in {\Z}.
\end{eqnarray}
Denote by $G_{l}(W_{1},W_{2})_{(h)}$ the space of all homogeneous
generalized left intertwining operators of weight $h$ from $W_{1}$ to
$W_{2}$ and set
\begin{eqnarray}
G^{o}_{l}(W_{1},W_{2})=\oplus_{h\in {\C}}G_{l}(W_{1},W_{2})_{(h)}.
\end{eqnarray}

Let  $E_{l}(W_{1},W_{2})$
be the space consisting of those
$\Phi(z)\in \left({\rm Hom}_{{\C}}(W_{1},W_{2})\right)\{z\}$
which satisfy conditions (Gl1) and (Gl3).

For any $a\in V, \Phi(z)\in E_{l}(W_{1},W_{2})$, we define               
\begin{eqnarray}\label{e2.27}
& &\!\!Y(a,z_{0})\circ_{l} \Phi(z_{2}):=\nonumber\\
&=&\!\!{\rm Res}_{z_{1}}\!\left(\!z_{0}^{-1}\delta\!\left(\!\frac{z_{1}-z_{2}}{z_{0}}\!
\right)Y_{W_{2}}(a,z_{1})\Phi (z_{2})-
z_{0}^{-1}\delta\!\left(\!\frac{z_{2}-z_{1}}{-z_{0}}\!\right)\Phi (z_{2})
Y_{W_{1}}(a,z_{1})\!\right)\!,
\end{eqnarray}
or equivalently
\begin{eqnarray}\label{e3.11}
a_{m}\circ_{l} \Phi (z_{2})={\rm Res}_{z_{1}}\!\left(\! (z_{1}-z_{2})^{m}
Y_{W_{2}}(a,z_{1})\Phi (z_{2})-(-z_{2}+z_{1})^{m}\Phi(z_{2})Y_{W_{1}}(a,z_{1})
\!\right)
\end{eqnarray}
for any $m\in {\Z}$. Then we have ([LL])

\bt{t3.2} With the action defined by (\ref{e2.27}),
the spaces $E_{l}(W_{1},W_{2})$ and $G_{l}(W_{1},W_{2})$ are
weak $V$-modules containing $G_{l}^{o}(W_{1},W_{2})$ as a weak submodule.
\et

Define a linear map $F_{l}(\cdot,z)$ from $G_{l}(W_{1},W_{2})$ to 
$\left({\rm Hom}_{{\C}}(W_{1},W_{2})\right)\{z\}$ as follows:
\begin{eqnarray}
F_{l}(\Phi,z)u_{1}=\Phi (z)u_{1}\;\;\;\mbox{ for }\Phi \in G_{l}(W_{1},W_{2}), u_{1}\in W_{1}.
\end{eqnarray}
For $a\in V, \Phi \in G_{l}(W_{1},W_{2})$, we have
\begin{eqnarray}
& &F_{l}(Y(a,z_{0})\circ_l\Phi,z_{2})\nonumber\\
&=&(Y(a,z_{0})\circ_l\Phi)(z)|_{z=z_{2}}\nonumber\\
&=&{\rm Res}_{z_{1}}\left(z_{0}^{-1}\delta\left(\frac{z_{1}-x}{z_{0}}
\right)Y(a,z_{1})\Phi(z)-z_{0}^{-1}\delta\left(\frac{z-z_{1}}
{-z_{0}}\right)\Phi(z)Y(a,z_{1})\right)|_{z=z_{2}}\nonumber\\
&=&{\rm Res}_{z_{1}}\left(z_{0}^{-1}\delta\!\left(\frac{z_{1}-z_{2}}{z_{0}}
\right)Y(a,z_{1})F_{l}(\Phi,z_{2})-z_{0}^{-1}\delta\!\left(\frac{z_{2}-z_{1}}
{-z_{0}}\right)F_{l}(\Phi,z_{2})Y(a,z_{1})\right).\nonumber\\
& &\mbox{}
\end{eqnarray}
It is well-known that this iterate formula implies associativity
[FHL]. Furthermore, (Gl3) gives the commutativity for
$F_{l}(\cdot,x)$. Therefore, $F(\cdot,z)$ satisfies the Jacobi
identity ([DL], [FHL], [L1]).  If $\Phi(x)\in G_{l}(W_{1},W_{2})$, by
definition we have
\begin{eqnarray}
F_{l}(L(-1)\circ_l\Phi,z)u_{1}=(L(-1)\circ_l \Phi(z))u_{1}={d\over
dz}\Phi(z)u_{1}={d\over dz}F_{l}(\Phi,z)u_{1}.\nonumber\\ 
\mbox{}
\end{eqnarray}
Therefore, $F_{l}(\cdot,z)$  is an intertwining 
operator of
type $\left(\!\begin{array}{c}W_{2}\\G_{l}(W_{1},W_{2}) W_{1}\end{array}\!\right)$.
It is clear that $F_{l}(\cdot,z)$ is injective in the sense that $F_{l}(\Phi,z)=0$ implies
$\Phi=0$ for $\Phi\in G_{l}(W_{1},W_{2})$. Then we have [LL]

\bt{t3.3} 
Let $W_{1}$ and $W_{2}$ be weak $V$-modules in the category $\cal{D}$. Then
for any weak $V$-module $M$ in the category $\cal{D}$ and any intertwining 
operator $I(\cdot,z)$
of type $\left(\!\begin{array}{c}W_{2}\\M W_{1}\end{array}\!\right)$,
there exists a unique $V$-homomorphism $f$ from $M$ to 
$G_{l}(W_{1},W_{2})$ such that $I(u,z)=F_{l}(f(u),x)$ for $u\in M$. Consequently,
the linear space ${\rm Hom}_{V}(M, G_{l}(W_{1},W_{2}))$ is linearly isomorphic to
$I\!\left(\!\begin{array}{c}W_{2}\\M W_{1} \end{array}\!\right)$.
\et

\begin{rem}\label{r3.4}
Suppose that $V$ satisfies the following conditions: (1) There are only
finitely
many inequivalent irreducible $V$-modules; (2) Any $V$-module is
completely reducible; (3) The fusion rules for any three modules are 
finite.  Then it was proved in [LL] that
for any pair of $V$-modules $W_{1}$ and $W_{2}$, 
there exists
a unique maximal module $\Delta(W_{1},W_{2})$ in $G_{l}(W_{1},W_{2})$, the contragredient 
module of which is isomorphic to
the tensor product of $W_{1}$ with the contragredient 
module of $W_2.$
\end{rem}

Let $I(\cdot,z)$ be an intertwining operator of type 
$\left(\!\begin{array}{c}W_{3}\\W_{1}W_{2}\end{array}\!\right)$. Then 
for any $u_{1}\in W_{1}$, $I(u_{1},z)$ is a generalized left intertwining 
operator from $W_{2}$ to $W_{3}$. Symmetrically, for any $u_{2}\in W_{2}$,
$I(\cdot,z)u_{2}$ is a generalized right intertwining operator from $W_{1}$
to $W_{3}$ as defined below.

\bd{d3.6} Let $W_{1}$ and $W_{2}$ be weak $V$-modules in the category ${\cal{D}}$.  
An element 
$\Phi (z)=\displaystyle{\sum_{\alpha\in {\C}}\Phi_{\alpha}z^{-\alpha-1}}
\in \left({\rm Hom}_{{\C}}(W_{1},W_{2})\right)\{z\}$ is called a 
{\em generalized right intertwining operator} from $W_{1}$ 
to $W_{2}$ if the following
conditions (Gr1)-(Gr3) hold: 

(Gr1)\hspace{0.25cm}  For any $\alpha\in {\C}, u_{1}\in W_{1}$,
$\Phi_{\alpha+n}u_{1}=0$ for $n\in {\Z}$ sufficiently large;

(Gr2)\hspace{0.25cm} $\Phi (z)L(-1)=\Phi'(z)\left(=\displaystyle{{d\over dz}}\Phi(z)\right);$

(Gr3)\hspace{0.25cm} For any $a\in V, u_{1}\in W_{1}$, there exists a positive integer
$k$ such that 
\begin{eqnarray}\label{erl}
(z_{0}+z_{2})^{k}Y_{W_{2}}(a,z_{0}+z_{2})\Phi(z_{2})u_{1}
=(z_{0}+z_{2})^{k}\Phi(z_{2})Y_{W_{1}}(a,z_{0})u_{1}.
\end{eqnarray}
\ed

Denote by $G_{r}(W_{1},W_{2})$ the space of all
generalized right intertwining operators from $W_{1}$ to
$W_{2}$. 
Homogeneous and primary generalized right intertwining operators can be defined similarly.
Denote by $G_{r}(W_{1},W_{2})_{(h)}$ the space of all weight-$h$
homogeneous generalized right intertwining operators from $W_{1}$ to
$W_{2}$ and set
\begin{eqnarray}
G^{0}_{r}(W_{1},W_{2})=\oplus_{h\in {\C}}G_{r}(W_{1},W_{2})_{(h)}.
\end{eqnarray}

For any $a\in V, \Phi(z)\in E_{r}(W_{1},W_{2})$, we define               
\begin{eqnarray}
& &\!\!Y(a,z_{1})\circ_{r} \Phi(z_{2}):=\nonumber\\
&=&\!\!Y_{W_{2}}(a,z_{1})\Phi (z_{2})-
{\rm Res}_{z_{0}}
z_{1}^{-1}\delta\left(\frac{z_{2}+z_{0}}{z_{1}}\right)\Phi (z_{2})
Y_{W_{1}}(a,z_{0}),
\end{eqnarray}
or equivalently
\begin{eqnarray}\label{er3.21}
a_{m}\circ_{r} \Phi (z_{2})=a_{m}\Phi (z_{2})
-{\rm Res}_{z_{0}}(z_{2}+z_{0})^{m}\Phi(z_{2})Y_{W_{1}}(a,z_{0})
\end{eqnarray}
for any $m\in {\Z}$. 

For any $\Phi(z)\in \left({\rm Hom}_{{\C}}(W_{1},W_{2})\right)\{z\}$, we define 
$T(\Phi(z))\in \left({\rm Hom}_{{\C}}(W_{1},W_{2})\right)\{z\}$
such that 
\begin{eqnarray}
T((\Phi(z))u_{1}=e^{zL(-1)}\Phi(e^{\pi i}z)u\;\;\;\mbox{ for }
u\in W_{1}.
\end{eqnarray}

\bp{p3.6}
Let $\Phi(z)\in \left({\rm Hom}_{{\C}}(W_{1},W_{2})\right)\{z\}$. Then 
$\Phi(z)\in E_{l}(W_{1},W_{2})$ if and only if
$T(\Phi(z))\in E_{r}(W_{1},W_{2})$.
\ep

{\bf Proof.} Since $T(\Phi(z))u=e^{zL(-1)}\Phi(e^{\pi i}z)u$ for any  
$u\in W_{1}$, we see that
$T(\Phi(z))$ satisfies (Gr1) if and only if $\Phi(z)$ satisfies (Gr1).
For $u\in W_{1}$, by definition we have
\begin{eqnarray}
& &{d\over dz}T(\Phi(z))u\nonumber\\
&=&L(-1)e^{zL(-1)}\Phi(e^{\pi i}z)+e^{zL(-1)}
{d\over dz}\Phi(e^{\pi i}z)u\nonumber\\
&=&e^{zL(-1)}[L(-1),\Phi(e^{\pi i}z)]+e^{zL(-1)}\Phi(e^{\pi i}z)L(-1)u+{d\over dz}
\Phi(e^{\pi i}z)u
\nonumber\\
&=&e^{zL(-1)}\left([L(-1),\Phi(e^{\pi i}z)]+{d\over dz}\Phi(e^{\pi i}z)u\right)
+T(\Phi(z))L(-1)u.
\end{eqnarray}
Then ${d\over dz}T(\Phi(z))=T(\Phi(z))L(-1)$ if and only if
$[L(-1),\Phi(e^{\pi i}z)]=-{d\over dz}\Phi(e^{\pi i}z)$. Thus $\Phi(z)$ satisfies (Gl2) 
if and only if
$T(\Phi(z))$ satisfies (Gr2).

For any $a\in V, u\in W_{1}$, suppose there is a positive integer $m$ such that
\begin{eqnarray}\label{e3.a1}
(z_{0}+z_{2})^{m}Y(a,z_{0}+z_{2})T(\Phi(z_{2}))u=(z_{0}+z_{2})^{m}T(\Phi(z_{2}))Y(a,z_{0})u.
\end{eqnarray}
Then
\begin{eqnarray}\label{e3.a2}
(z_{0}+z_{2})^{m}e^{-z_{2}L(-1)}Y(a,z_{0}+z_{2})T(\Phi(z_{2}))u
=(z_{0}+z_{2})^{m}e^{-z_{2}L(-1)}T(\Phi(z_{2}))Y(a,z_{0})u.
\end{eqnarray}
Thus
\begin{eqnarray}\label{e3.a3}
(z_{0}+z_{2})^{m}Y(a,z_{0})e^{-z_{2}L(-1)}T(\Phi(z_{2}))u
=(z_{0}+z_{2})^{m}e^{-z_{2}L(-1)}T(\Phi(z_{2}))Y(a,z_{0})u.
\end{eqnarray}
That is,
\begin{eqnarray}\label{e3.a4}
(z_{0}+z_{2})^{m}Y(a,z_{0})\Phi(e^{\pi i}z_{2}))u
=(z_{0}+z_{2})^{m}\Phi(e^{\pi i}z_{2}))Y(a,z_{0})u.
\end{eqnarray}
Therefore $\Phi(z)$ satisfies (Gl3). Since every step in the proof can
be reversed,
 $\Phi(z)$ satisfies (Gl3) if and only if $T(\Phi(z))$ satisfies (Gr3).
 $\;\;\;\;\Box$

\bp{p3.7} For any $a\in V, \Phi(z)\in E_{l}(W_{1},W_{2})$ we have
$$T(Y(a,z_{0})\circ_{l}\Phi(z))=Y(a,z_{0})\circ _{r}\Phi(z).$$
\ep

{\bf Proof.} By definition we have
\begin{eqnarray}
& &\;\;\;T(a_m\circ_{l} \Phi (z_{2}))\nonumber\\
& &={\rm Res}_{z_{1}}\!\left(\! (z_{1}+z_{2})^{m}e^{z_{2}L(-1)}
Y_{W_{2}}(a,z_{1})\Phi (e^{\pi i}z_{2})-(z_{2}+z_{1})^{m}e^{z_{2}L(-1)}\Phi(e^{\pi i}z_{2})
Y_{W_{1}}(a,z_{1})\!\right)\nonumber\\
& &={\rm Res}_{z_{1}}(z_{1}+z_{2})^{m}Y_{W_{2}}(a,z_{1}+z_{2})e^{z_{2}L(-1)}
\Phi (e^{\pi i}z_{2})\nonumber\\
& &\;\;\;-{\rm Res}_{z_{1}}(z_{2}+z_{1})^{m}e^{z_{2}L(-1)}\Phi(e^{\pi i}z_{2})Y_{W_{1}}(a,z_{1})\nonumber\\
& &={\rm Res}_{z_{1}}\!\left(\! (z_{1}+z_{2})^{m}Y_{W_{2}}(a,z_{1}+z_{2})
T(\Phi (z_{2}))-(z_{2}+z_{1})^{m}T(\Phi(z_{2}))Y_{W_{1}}(a,z_{1})
\!\right)\nonumber\\
& &=a_{m}T(\Phi(z_{2}))-{\rm Res}_{z_{1}}(z_{2}+z_{1})^{m}T(\Phi(z_{2}))Y_{W_{1}}(a,z_{1})
\nonumber\\
& &=a_m\circ_{r}\Phi(z_{2}).\;\;\;\;\Box
\end{eqnarray}

Define a linear map $T$ as follows:
\begin{eqnarray}
T: & &E_{l}(W_{1}, W_{2})\rightarrow E_{r}(W_{1}, W_{2});\nonumber\\
& &\Phi(z)\mapsto T(\Phi(z))\;\;\;\mbox{ for }\Phi(z)\in E_{l}(W_{1}, W_{2}).
\end{eqnarray}
Then we obtain

\bp{p3.8} The spaces $E_{r}(W_{1},W_{2})$ and $G_{r}(W_{1},W_{2})$ are
weak $V$-modules and the linear map $T$ is a $V$-isomorphism.
\ep

Similarly we have

\bt{t3.9} 
Let $W_{1}$and $W_{2}$ be weak $V$-modules in the category ${\cal{D}}$. Then
for any weak $V$-module $M$ in the category $\cal{D}$ and any intertwining 
operator $I(\cdot,z)$
of type $\left(\!\begin{array}{c}W_{2}\\W_{1}M\end{array}\!\right)$,
there exists a unique $V$-homomorphism $f$ from $M$ to 
$G_{r}(W_{1},W_{2})$ such that $I(u_{1},z)u=F_{r}(u_{1},z)f(u)$ for $u_{1}\in W_{1},u\in M$.
Consequently,
the linear space ${\rm Hom}_{V}(M, G_{r}(W_{1},W_{2}))$ is linearly isomorphic to
$I\!\left(\!\begin{array}{c}W_{2}\\W_{1} M\end{array}\!\right)$.
\et

\section{The associativity of tensor products}
The main goal of this section is to prove the associativity of tensor
products of three modules for a vertex operator algebra $V$ under
Huang's convergence condition.  Throughout this section, $V$ is a
fixed vertex operator algebra and we assume the existence of tensor
products in either of the categories ${\cal{C}}(V)$ or
${\cal{C}}^{1}(V)$. We denote either one by ${\cal D}.$

Propositions \ref{pu1} and \ref{pu2} immediately follow
from the universal property of tensor product in Definition \ref{d2.1}.

\bp{pu1} 
Let $W_{i}$ $(i=1,2,3)$ be weak $V$-modules in the category $\cal{D}$. 
Let $(W_{23},F_{23})$
be a tensor product of $(W_{2},W_{3})$ and let 
$(W_{1,23}, F_{1,23})$
be a tensor product of $(W_{1}, W_{23})$. Then
the quadruple $(W_{23}, W_{1,23}, F_{23}, F_{1,23})$ 
satisfies the following universal property: For any quadruple 
$(M_{1},M_{2}, I_{1}, I_{2})$, 
where 
$M_{1},M_{2}$ are weak $V$-modules in the category $\cal{D}$ and
$I_{1,23},\; I_{23}$ are any intertwining operators 
of types $\left(\!\begin{array}{c}M_{2}\\W_{1} M_{1}\end{array}\!\right)$,
$\left(\!\begin{array}{c}M_{1}\\W_{2}W_{3}\end{array}\!\right)$, respectively,
there exists a unique $V$-homomorphism $f$ from
 $W_{1,23}$ to $M_{2}$ such that
\begin{eqnarray}
fF_{1,23}(u_{1},z_{1})F_{23}(u_{2},z_{2})u_{3}
=I_{1,23}(u_{1},z_{1})I_{23}(u_{2},z_{2})u_{3}
\end{eqnarray}
for $u_{i}\in W_{i}$ $(i=1,2,3)$. 
\ep

It is easy to see that this universal property characterizes $W_{1,23}$ uniquely.
Similarly, we have

\bp{pu2}
Let $(W_{12}, F_{12})$ be a tensor product of 
$(W_{1},W_{2})$ and let $(W_{12,3}, F_{12,3})$ 
is a tensor product of $(W_{12}, W_{3})$. Then the quadruple
$(W_{12}, W_{12,3}, F_{12}, 
F_{12,3})$ satisfies the following  corresponding universal property: For any
quadruple $(M_{1},M_{2}, I_{12}, I_{12,3})$, where $M_{1},\;M_{2}$ are 
weak $V$-modules in the category $\cal{D}$ and $I_{12},\; I_{12,3}$ are 
intertwining operators of types
$\left(\!\begin{array}{c}M_{1}\\W_{1}W_{2}\end{array}\!\right)$ and
$\left(\!\begin{array}{c}M_{2}\\M_{1}W_{3}\end{array}\!\right)$, respectively,
there exists a unique $V$-homomorphism $g$ from 
$W_{12,3}$ to $M_{2}$ such that
\begin{eqnarray}
gF_{12,3}(F_{12}(u_{1},z_{0})u_{2},z_{2})u_{3}=
I_{12,3}(I_{12}(u_{1},z_{0})u_{2},z_{2})w_{3}
\end{eqnarray}
for any $u_{i}\in W_{i}\; (i=1,2,3)$.
\ep

Our main goal is to prove that the tensor product modules $W_{1,23}$ 
and $W_{12,3}$ are isomorphic under Huang's convergence assumptions [H1]-[H2]:

{\bf Convergence and extension property for products:} Let $W_{i}$ $(i=1,2,3,4)$ be weak $V$-modules in the category $\cal{D}$ and 
let ${\cal{Y}}_{1}$, 
${\cal{Y}}_{2}$ be intertwining operators of types
$\left(\!\begin{array}{c}W_{4}\\W_{1}W_{5}\end{array}\!\right)$ and 
$\left(\!\begin{array}{c}W_{5}\\W_{2}W_{3}\end{array}\!\right)$, respectively.
There exists an integer $N$ such that
for homogeneous $w_{i}\in W_{i}\; (i=1,2,3), w_{4}'\in W_{4}'$, there exist 
$k\in {\N}, r_{i}, s_{i}\in {\C}\;(i=1,2,\cdots,k)$ and analytic functions 
$f_{i}(z)$ on $|z|<1$ such that
${\rm Re} ({\rm gwt}w_{1}+{\rm gwt}w_{2}+s_{i})>N$ and that
$$\<w_{4}',{\cal{Y}}_{1}(w_{1},z_{1}){\cal{Y}}_{2}(w_{2},z_{2})w_{3}\>$$
converges in the domain $|z_{1}|>|z_{2}|>0$ to a function which can be analytically
extended to the multi-valued function
\begin{eqnarray}\label{e3.6}
\sum_{i=1}^{k}z_{2}^{r_{i}}(z_{1}-z_{2})^{s_{i}}f_{i}\left(\frac{z_{1}-z_{2}}{z_{2}}\right)
\end{eqnarray}
in the domain $|z_{2}|>|z_{1}-z_{2}|>0$. 

Let $f(z_{1},z_{2})=\sum_{m,n}a_{m,n}z_{1}^{m}z_{2}^{n}$ be a formal series 
which is convergent in the domain $|z_{1}|>|z_{2}|>0$ and can be analytically extended to
the multi-valued function (\ref{e3.6}) in the domain $|z_{2}|>|z_{1}-z_{2}|>0$. 
Then we define
\begin{eqnarray}
E_{12}f(z_{1},z_{2})=\sum_{i=1}^{k}z_{2}^{r_{i}}(z_{1}-z_{2})^{s_{i}}f_{i}
\left(\frac{z_{1}-z_{2}}{z_{2}}\right)
\end{eqnarray}
as a formal series. It is clear that 
$$E_{12}P(z_{1},z_{2})f(z_{1},z_{2})=P(z_{1},z_{2})E_{12}f(z_{1},z_{2})
\;\;\;\mbox{for any }P(z_{1},z_{2})\in \C[z_{1}^{\pm},z_{2}^{\pm}].$$
That is, $E_{12}$ is a $\C[z_{1}^{\pm},z_{2}^{\pm}]$-homomorphism.

{\bf Convergence and extension property for iterates:} Let ${\cal{Y}}_{3}$ and ${\cal{Y}}_{4}$ be intertwining operators of types
$\left(\!\begin{array}{c}W_{5}\\W_{1}W_{2}\end{array}\!\right)$ and 
$\left(\!\begin{array}{c}W_{4}\\W_{5}W_{3}\end{array}\!\right)$, respectively. Then
there exists an integer $\tilde{N}$ depending only ${\cal{Y}}_{3}$ and ${\cal{Y}}_{4}$, and 
for homogeneous $w_{i}\in W_{i}$ $(i=1,2,3),$  $w_{4}'\in W_{4}'$, there exist 
$k\in {\N}, \tilde{r}_{i}, \tilde{s}_{i}\in {\C}, i=1,\cdots,k$ and analytic functions 
$\tilde{f}_{i}(z)$ on $|z|<1, i=1,\cdots,k$, satisfying
\begin{eqnarray}
{\rm Re}({\rm gwt}w_{2}+{\rm gwt}w_{3}+\tilde{s}_{i})>\tilde{N},\;\; i=1,\cdots, k,
\end{eqnarray}
such that
\begin{eqnarray}
\< w_{4}',{\cal{Y}}_{4}({\cal{Y}}_{3}(w_{1},z_{1}-z_{2})w_{2},z_{2})w_{3}\>
\end{eqnarray}
is convergent when $|z_{2}|>|z_{1}-z_{2}|>0$ and can be analytically extended to the 
multi-valued function
\begin{eqnarray}
\sum_{i=1}^{k}z_{1}^{\tilde{r}_{i}}z_{2}^{\tilde{s}_{i}}\tilde{f}_{i}\left({z_{2}\over z_{1}}
\right)\end{eqnarray}
when $|z_{1}|>|z_{2}|>0$. For convenience we write
\begin{eqnarray}
E^{12}\< w_{4}',{\cal{Y}}_{4}({\cal{Y}}_{3}(w_{1},z_{1}-z_{2})w_{2},z_{2})w_{3}\>
=\sum_{i=1}^{k}z_{1}^{\tilde{r}_{i}}z_{2}^{\tilde{s}_{i}}\tilde{f}_{i}
\left({z_{2}\over z_{1}}\right)
\end{eqnarray}
as formal power series. The following proposition can be found in [H1].

\bp{ph}
The convergence and extension property for products is equivalent to the 
convergence and extension property for iterates.
\ep

\begin{rem}\label{r3.9}
For most of the known examples, one can check that the convergence
follows from the KZ equations ([KZ], [TK]) or from free field realizations
(cf. [BF], [Fe]).
\end{rem}

For any fixed $w_{i}\in W_{i}\; (i=1,2)$, we define a linear map
$\psi(w_{1},w_{2};z_{0},z_{2})$ from $W_{3}$ to $W_{4}\{z_{0},z_{2}\}$ as follows:
\begin{eqnarray}
& &\;\;\;\< w_{4}',\psi(w_{1},w_{2};z_{0},z_{2})w_{3}\>:=\nonumber\\
& &=E_{12}\< w_{4}', {\cal{Y}}_{1}(w_{1},z_{1})
{\cal{Y}}_{2}(w_{2},z_{2})w_{3}\>|_{z_{1}=z_{0}+z_{2}}
\nonumber\\
& &=\sum_{i=1}^{k}z_{2}^{r_{i}}z_{0}^{s_{i}}f_{i}\left({z_{0}\over z_{2}}\right).
\end{eqnarray}
By definition, $\< w_{4}',\psi(w_{1},w_{2};z_{0},z_{2})w_{3}\>$ is a formal series 
in $z_{0}$ and $z_{2}$, which is convergent in the domain 
$|z_{2}|>|z_{0}|>0$.

For any $w_{1}\in W_{1},w_{2}\in W_{2}$, we set
\begin{eqnarray}
\psi(w_{1},w_{2};z_{0},z_{2})=\sum_{n\in {\C}}\psi_{n}(w_{1},w_{2},z_{2})z_{0}^{-n-1}
=\sum_{m,n\in {\C}}\psi_{n}(w_{1},w_{2})_{m}z_{0}^{-n-1}z_{2}^{-m-1}.
\end{eqnarray}

\begin{rem}\label{r3.10}
{} From the assumption that ${\rm Re}({\rm gwt}w_{1}+{\rm
gwt}w_{2}+s_{i})>N$ for $i=1,\cdots,k$, we see that there exists a
real number $h$ such that $\psi_{n}(w_{1},w_{2},z_{2})=0$ whenever
${\rm Re}\,n>h$.
\end{rem}

\bp{p3.11}
For any $w_{1}\in W_{1},w_{2}\in W_{2}, n\in {\C}$, we have
$$\psi_{n}(w_{1},w_{2},z_{2})\in G_{l}^{o}(W_{3},W_{4}).$$
\ep

{\bf Proof.} For homogeneous vectors $w_{i}\in W_{i}\;(i=1,2,3), w_{4}'\in W_{4}'$, we have
\begin{eqnarray}
& &\;\;\;\<w_{4}',[L(0),\psi(w_{1},w_{2};z_{0},z_{2})]w_{3}\>\nonumber\\
& &=\<w_{4}',L(0)\psi(w_{1},w_{2};z_{0},z_{2})w_{3}-\psi(w_{1},w_{2};z_{0},z_{2})L(0)w_{3}\>
\nonumber\\
& &=\<L(0)w_{4}',\psi(w_{1},w_{2};z_{0},z_{2})w_{3}\>-\<w_{4}',\psi(w_{1},w_{2};z_{0},z_{2})L(0)w_{3}\>
\nonumber\\
& &=E_{12}
\<L(0)w_{4}',{\cal{Y}}_{1}(w_{1},z_{1}){\cal{Y}}_{2}(w_{2},z_{2})w_{3}\>|_{z_{1}=z_{0}+z_{2}}
\nonumber\\
& &\;\;\; -E_{12}\<w_{4}',{\cal{Y}}_{1}(w_{1},z_{1}){\cal{Y}}_{2}(w_{2},z_{2})
L(0)w_{3}\>|_{z_{1}=z_{0}+z_{2}}
\nonumber\\
& &=E_{12}
\<w_{4}',L(0){\cal{Y}}_{1}(w_{1},z_{1}){\cal{Y}}_{2}(w_{2},z_{2})w_{3}\>|_{z_{1}=z_{0}+z_{2}}
\nonumber\\
& &\;\;\; -E_{12}\<w_{4}',{\cal{Y}}_{1}(w_{1},z_{1}){\cal{Y}}_{2}(w_{2},z_{2})
L(0)w_{3}\>|_{z_{1}=z_{0}+z_{2}}\nonumber\\
& &=E_{12}
\<w_{4}',z_{1}{\partial\over \partial z_{1}}
{\cal{Y}}_{1}(w_{1},z_{1}){\cal{Y}}_{2}(w_{2},z_{2})w_{3}\>|_{z_{1}=z_{0}+z_{2}}
\nonumber\\
& &\;\;+E_{12}\<w_{4}',{\cal{Y}}_{1}(L(0)w_{1},z_{1}){\cal{Y}}_{2}(w_{2},z_{2})
w_{3}\>|_{z_{1}=z_{0}+z_{2}}
\nonumber\\
& &\;\;+E_{12}
\<w_{4}',z_{2}{\partial\over \partial z_{2}}
{\cal{Y}}_{1}(w_{1},z_{1}){\cal{Y}}_{2}(w_{2},z_{2})w_{3}\>|_{z_{1}=z_{0}+z_{2}}
\nonumber\\
& &\;\;+E_{12}\<w_{4}',{\cal{Y}}_{1}(w_{1},z_{1}){\cal{Y}}_{2}(L(0)w_{2},z_{2})
w_{3}\>|_{z_{1}=z_{0}+z_{2}}
\nonumber\\
& &=E_{12}
\<w_{4}',(z_{0}+z_{2}){\partial\over \partial z_{1}}
{\cal{Y}}_{1}(w_{1},z_{1}){\cal{Y}}_{2}(w_{2},z_{2})w_{3}\>|_{z_{1}=z_{0}+z_{2}}
\nonumber\\
& &\;\;+E_{12}
\<w_{4}',z_{2}{\partial\over \partial z_{2}}
{\cal{Y}}_{1}(w_{1},z_{1}){\cal{Y}}_{2}(w_{2},z_{2})w_{3}\>|_{z_{1}=z_{0}+z_{2}}
\nonumber\\
& &\;\;+\<w_{4}',\psi(L(0)w_{1},w_{2};z_{0},z_{2})w_{3}+\psi(w_{1},L(0)w_{2};z_{0},z_{2})w_{3}\>
\nonumber\\
& &=z_{0}{\partial\over \partial z_{0}}\<w_{4}',\psi(w_{1},w_{2};z_{0},z_{2})w_{3}\>\nonumber\\
& &\;\;+z_{2}{\partial\over \partial z_{2}}\<w_{4}',\psi(w_{1},w_{2};z_{0},z_{2})w_{3}\>
\nonumber\\
& &\;\;+\<w_{4}',\psi(L(0)w_{1},w_{2};z_{0},z_{2})w_{3}+\psi(w_{1},L(0)w_{2};z_{0},z_{2})w_{3}\>.
\end{eqnarray}
Thus
\begin{eqnarray}
& &\ \ \ [L(0), \psi_{n}(w_{1},w_{2})_{m}]\nonumber\\
& &=(-n-1-m-1)\psi_{n}(w_{1},w_{2})_{m}+
\psi_{n}(L(0)w_{1},w_{2})_{m}+\psi_{n}(w_{1},L(0)w_{2})_{m}.
\end{eqnarray}
Suppose that $w_{1}, w_{2}, w_{3}$ have generalized weights $h_{1},h_{2}, h_{3}$, respectively. 
Then we obtain
\begin{eqnarray}
&&\ \ \ (L(0)-h_{1}-h_{2}-h_{3}+n+1+m+1)\psi_{n}(w_{1},w_{2})_{m}w_{3}\nonumber\\
& &=\psi_{n}((L(0)-h_{1})w_{1},w_{2})_{m}w_{3}
+\psi_{n}(w_{1},(L(0)-h_{2})w_{2})_{m}w_{3}\nonumber\\
& &+\psi_{n}(w_{1},w_{2})_{m}(L(0)-h_{3})w_{3}.
\end{eqnarray}
Then it follows that $\psi_{n}(w_{1},w_{2})_{m}$ is homogeneous of
generalized weight $h_{1}+h_{2}-n-1-m-1$ for $m,n\in \C$. (Note that
$\psi_{n}(w_{1},w_{2})_{m}$ is an eigenvector for $L(0)$ if both
$w_1$ and $w_2$ are.) Since the
generalized weights of $W_{4}$ are truncated from below by definition
(recall that ${\cal D}={\cal C}(V)$ or ${\cal C}^1(V)$), it follows
that $\psi_{n}(w_{1},w_{2},z_{2})$ satisfies (Gl1).  Similarly, we
have
\begin{eqnarray}
& &\;\;\;\<w_{4}',[L(-1), \psi(w_{1},w_{2};z_{0},z_{2})]w_{3}\>\nonumber\\
& &=E_{12}
\<w_{4}',\left({\partial\over \partial z_{1}}+{\partial\over \partial z_{2}}\right)
{\cal{Y}}_{1}(w_{1},z_{1}){\cal{Y}}_{2}(w_{2},z_{2})w_{3}\>|_{z_{1}=z_{0}+z_{2}}
\nonumber\\
& &={\partial\over \partial z_{2}}E_{12}
\<w_{4}',{\cal{Y}}_{1}(w_{1},z_{1}){\cal{Y}}_{2}(w_{2},z_{2})w_{3}\>|_{z_{1}=z_{0}+z_{2}}
\nonumber\\
& &={\partial\over \partial z_{2}}\<w_{4}',\psi(w_{1},w_{2};z_{0},z_{2})w_{3}\>.
\end{eqnarray}
Then $\psi_{n}(w_{1},w_{2},z_{2})$ satisfies (Gl2).
Suppose that $\psi_{m}(w_{1},w_{2},z_{2})=0$ whenever Re\,$m>h$ for some real number $h$. Then
$${\rm Res}_{z_{0}}z_{0}^{m}f(z_{0},z_{2})\psi(w_{1},w_{2};z_{0},z_{2})=0$$
for any $f(z_{0},z_{2})\in {\C}[z_{0},z_{2}]$. For any $a\in V, n\in {\C}$, 
suppose that $k$ is a positive integer such that
$k+n>h$ and that
\begin{eqnarray}
& &(z_{1}-z_{2})^{k}Y(a,z_{1}){\cal{Y}}_{1}(w_{1},z_{2})
=(z_{1}-z_{2})^{k}{\cal{Y}}_{1}(w_{1},z_{2})Y(a,z_{1}),\\
& &(z_{1}-z_{2})^{k}Y(a,z_{1}){\cal{Y}}_{2}(w_{1},z_{2})
=(z_{1}-z_{2})^{k}{\cal{Y}}_{2}(w_{1},z_{2})Y(a,z_{1}).
\end{eqnarray}
Then we have
\begin{eqnarray}
& &\;\;\;\;(z-z_{2})^{3k}\<w_{4}', Y(a,z)\psi_{n}(w_{1},w_{2},z_{2})w_{3}\>
\nonumber\\
& &={\rm Res}_{z_{0}}z_{0}^{n}(z-z_{2})^{3k}
 \<w_{4}', Y(a,z)\psi(w_{1},w_{2};z_{0},z_{2})w_{3}\>\nonumber\\
& &=\sum_{i=0}^{2k}{\rm Res}_{z_{0}}{2k\choose i}(z-z_{0}-z_{2})^{i}
z_{0}^{n+2k-i}(z-z_{2})^{k}
\<w_{4}',Y(a,z)\psi(w_{1},w_{2};z_{0},z_{2})w_{3}\>\nonumber\\
& &=\sum_{i=k}^{2k}{\rm Res}_{z_{0}}{2k\choose i}(z-z_{0}-z_{2})^{i}
z_{0}^{n+2k-i}(z-z_{2})^{k}\<w_{4}',Y(a,z)\psi(w_{1},w_{2};z_{0},z_{2})w_{3}\>\nonumber\\
& &=\sum_{i=k}^{2k}{\rm Res}_{z_{0}}{2k\choose i}(z-z_{0}-z_{2})^{i}
z_{0}^{n+2k-i}(z-z_{2})^{k}\cdot \nonumber\\
& &\;\;\;\cdot E_{12}
\<w_{4}', Y(a,z){\cal{Y}}_{1}(w_{1},z_{1}){\cal{Y}}_{2}(w_{2},z_{2})w_{3}\>|_{z_{1}=z_{0}+z_{2}}
\nonumber\\
& &=\sum_{i=k}^{2k}{\rm Res}_{z_{0}}{2k\choose i}\cdot \nonumber\\
& &\;\;\;\cdot E_{12}
z_{0}^{n+2k-i}
\<w_{4}', (z-z_{1})^{i}(z-z_{2})^{k}{\cal{Y}}_{1}(w_{1},z_{1}){\cal{Y}}_{2}(w_{2},z_{2})
Y(a,z)w_{3}\>|_{z_{1}=z_{0}+z_{2}}
\nonumber\\
& &=\sum_{i=k}^{2k}{\rm Res}_{z_{0}}{2k\choose i}(z-z_{0}-z_{2})^{i}
z_{0}^{n+2k-i}(z-z_{2})^{k}\<w_{4}', 
\psi(w_{1},w_{2};z_{0},z_{2})Y(a,z)w_{3}\>\nonumber\\
& &=\sum_{i=0}^{2k}{\rm Res}_{z_{0}}{2k\choose i}(z-z_{0}-z_{2})^{i}
z_{0}^{n+2k-i}(z-z_{2})^{k}\<w_{4}',\psi(w_{1},w_{2};z_{0},z_{2})Y(a,z)w_{3}\>\nonumber\\
& &={\rm Res}_{z_{0}}z_{0}^{n}(z-z_{2})^{3k}
\<w_{4}', \psi(w_{1},w_{2};z_{0},z_{2})Y(a,z)w_{3}\>\nonumber\\
& &=(z-z_{2})^{3k}\<w_{4}',\psi_{n}(w_{1},w_{2},z_{2})Y(a,z)w_{3}\>.
\end{eqnarray}
Thus $\psi_{n}(w_{1},w_{2},z_{2})\in G^{0}_{l}(W_{3},W_{4})$ for 
$n\in {\C}, w_{1}\in W_{1},w_{2}\in W_{2}$. 
$\;\;\;\;\Box$

Let $W$ be the subspace of $ G^{0}_{l}(W_{3},W_{4})$ linearly spanned by all 
$\psi_{n}(w_{1},w_{2},z_{2})$ for $n\in {\C}, w_{1}\in W_{1},w_{2}\in W_{2}$.

\begin{lem}\label{l3.2}
For any $a\in V, w_{1}\in W_{1},w_{2}\in W_{2}$, we have
\begin{eqnarray}
a_{m}\circ_{l} \psi(w_{1},w_{2};z_{0},z_{2})
=\sum_{i=0}^{\infty}{m\choose i}\psi(a_{i}w_{1},w_{2};z_{0},z_{2})z_{0}^{m-i}
-\psi(w_{1},a_{m}w_{2};z_{0},z_{2}).
\end{eqnarray}
\end{lem}

{\bf Proof.} For any $a\in V, w_{i}\in W_{i}\;(i=1,2,3), w_{4}'\in W_{4}'$, 
we have
\begin{eqnarray}
& &\;\;\;\;\<w_{4}',a_{m}\circ_{l} \psi(w_{1},w_{2};z_{0},z_{2})w_{3}\>
\nonumber\\
& &={\rm Res}_{z_{1}}\<w_{4}',(z_{1}-z_{2})^{m}Y(a,z_{1})
\psi(w_{1},w_{2};z_{0},z_{2})w_{3}\>\nonumber\\
& &\;\;\;-{\rm Res}_{z_{1}}\<w_{4}',(-z_{2}+z_{1})^{m}
\psi(w_{1},w_{2};z_{0},z_{2})Y(a,z_{1})w_{3}\>\nonumber\\
& &={\rm Res}_{z_{1}}E_{32}
\<w_{4}',(z_{1}-z_{2})^{m}
Y(a,z_{1}){\cal{Y}}_{1}(w_{1},z_{3}){\cal{Y}}_{2}(w_{2},z_{2})w_{3}\>|_{z_{3}=z_{0}+z_{2}}
\nonumber\\
& &\;\;\;-{\rm Res}_{z_{1}}E_{32}
\<w_{4}', (-z_{2}+z_{1})^{m}{\cal{Y}}_{1}(w_{1},z_{3})
{\cal{Y}}_{2}(w_{2},z_{2})Y(a,z_{1})w_{3}\>|_{z_{3}=z_{0}+z_{2}}\nonumber\\
& &={\rm Res}_{z_{1}}{\rm Res}_{z_{4}}(z_{1}-z_{2})^{m}
z_{1}^{-1}\delta\left(\frac{z_{0}+z_{2}+z_{4}}{z_{1}}\right)\nonumber\\
& &\;\;\;\;\cdot E_{32}\<w_{4}', 
{\cal{Y}}_{1}(Y(a,z_{4})w_{1},z_{3}){\cal{Y}}_{2}(w_{2},z_{2})w_{3}\>|_{z_{3}=z_{0}+z_{2}}
\nonumber\\
& &\;\;\;-E_{32}\<w_{4}',{\cal{Y}}_{1}(w_{1},z_{4})
{\cal{Y}}_{2}(a_{m}w_{2},z_{2})w_{3}\>|_{z_{3}=z_{0}+z_{2}}\nonumber\\
& &={\rm Res}_{z_{4}}\iota_{04}E_{32}
\<w_{4}',(z_{0}+z_{4})^{m}{\cal{Y}}_{1}(Y(a,z_{4})w_{1},z_{3})
{\cal{Y}}_{2}(w_{2},z_{2})w_{3}\>|_{z_{3}=z_{0}+z_{2}}\nonumber\\
& &\;\;\;-\<w_{4}',\psi(w_{1},a_{m}w_{2};z_{0},z_{2})w_{3}\>\nonumber\\
& &=\sum_{i=0}^{\infty}{m\choose i}\<w_{4}',\psi(a_{i}w_{1},w_{2};z_{0},z_{2})
w_{3}\>z_{0}^{m-i}
-\<w_{4}',\psi(w_{1},a_{m}w_{2};z_{0},z_{2})w_{3}\>.
\end{eqnarray}
Now the lemma follows immediately.$\;\;\;\;\Box$

It follows from Lemma \ref{l3.2} that $W$ is a weak $V$-submodule of 
$G_{l}^{o}(W_{1},W_{2})$. Then we 
obtain a natural intertwining 
operator
${\cal{Y}}_{5}(\cdot,z)$ of type
$\left(\!\begin{array}{c}W_{4}\\W W_{3}\end{array}\!\right)$ such that
${\cal{Y}}_{5}(w,z)=F_{l}(w,z)$, where $F_{l}(\cdot,z)$ is the natural intertwining 
operator of type
$\left(\!\begin{array}{c}W_{4}\\G_{l}^{o}(W_{3},W_{4})W_{3}\end{array}\!\right)$.

For any $w_{1}\in W_{1}, w_{2}\in W_{2}$,
we define ${\cal{Y}}_{4}(w_{1},z_{0})w_{2}\in W\{z_{0}\}$ such that
\begin{eqnarray}
\left({\cal{Y}}_{4}(w_{1},z_{0})w_{2}\right)(z_{2})=\psi(w_{1},w_{2};z_{0},z_{2})
=\sum_{n\in{\C}}\psi_{n}(w_{1},w_{2},z_{2})z_{0}^{-n-1}. 
\end{eqnarray}

\bp{p3.3x}
The linear map ${\cal{Y}}_{4}$ is an intertwining operator of type
$\!\left(\!\!\begin{array}{c}W\\W_{1} W_{2}\end{array}\!\!\right)$.
\ep

{\bf Proof.} For any $w_{1}\in W_{1},w_{2}\in W_{2}$, it follows from
the assumption that there exists a real number $h$ such that
$\psi_{n}(w_{1},w_{2},z_{2})=0$ whenever ${\rm Re}\,n>h$.  Then it
follows from the definition that ${\cal{Y}}_{4}$ satisfies the
truncation condition (I1).  For $w_{i}\in W_{i}\;(i=1,2,3), w_{4}'\in
W_{4}'$, by definition we again have
\begin{eqnarray}
& &\;\;\;\<w_{4}',\left({\cal{Y}}_{4}(L(-1)w_{1},z_{0})w_{2}\right)(z_{2})w_{3}\>\nonumber\\
& &=\<w_{4}',\psi(L(-1)w_{1},w_{2};z_{0},z_{2})w_{3}\>\nonumber\\
& &=E_{12}
\<w_{4}',{\cal{Y}}_{1}(L(-1)w_{1},z_{1}){\cal{Y}}_{2}(w_{2},z_{2})w_{3}\>|_{z_{1}=z_{0}+z_{2}}
\nonumber\\
& &=E_{12}
\<w_{4}',\left({\partial\over \partial z_{1}} {\cal{Y}}_{1}(w_{1},z_{1})
{\cal{Y}}_{2}(w_{2},z_{2})\right)w_{3}\>|_{z_{1}=z_{0}+z_{2}}\nonumber\\
& &={\partial\over \partial z_{0}}E_{12}
\<w_{4}', {\cal{Y}}_{1}(w_{1},z_{1})
{\cal{Y}}_{2}(w_{2},z_{2})w_{3}\>|_{z_{1}=z_{0}+z_{2}}\nonumber\\
& &={\partial\over \partial z_{0}}\<w_{4}',\left({\cal{Y}}_{4}(w_{1},z_{0})w_{2}\right)(z_{2})w_{3}\>.
\end{eqnarray}
Then ${\cal{Y}}_{4}(\cdot,z)$ satisfies the $L(-1)$-derivative property (I2). It remains to 
prove the Jacobi identity.
For $a\in V, w_{i}\in W_{i}$ $(i=1,2,3), w_{'}\in W_{4}'$, we have
\begin{eqnarray}\label{e3.21}
& &\;\;\;\<w_{4}',\left(Y(a,z){\cal{Y}}_{4}(w_{1},z_{0})w_{2}\right)(z_{2})w_{3}\>\nonumber\\
& &=\<w_{4}',\left(Y(a,z)\circ_{l} \psi(w_{1},w_{2};z_{0},z_{2})\right)w_{3}\>\nonumber\\
& &=\<w_{4}',{\rm Res}_{z_{1}} z^{-1}\delta\left(\frac{z_{1}-z_{2}}{z}\right)
Y(a,z_{1})\psi(w_{1},w_{2};z_{0},z_{2})w_{3}\>\nonumber\\
& &\;\;\;-\<w_{4}', {\rm Res}_{z_{1}}z^{-1}\delta\left(\frac{z_{2}-z_{1}}{-z}\right)
\psi(w_{1},w_{2};z_{0},z_{2})Y(a,z_{1})w_{3}\>\nonumber\\
& &=E_{32}
\<w_{4}',{\rm Res}_{z_{1}} z^{-1}\delta\left(\frac{z_{1}-z_{2}}{z}\right)
Y(a,z_{1}){\cal{Y}}_{1}(w_{1},z_{3}){\cal{Y}}_{2}(w_{2},z_{2})w_{3}\>|_{z_{3}=z_{0}+z_{2}}\nonumber\\
& &\;\;\;-E_{32}
\<w_{4}',{\rm Res}_{z_{1}}z^{-1}\delta\left(\frac{z_{2}-z_{1}}{-z}\right)
{\cal{Y}}_{1}(w_{1},z_{3}){\cal{Y}}_{2}(w_{2},z_{2})Y(a,z_{1})w_{3}\>|_{z_{3}=z_{0}+z_{2}},\;\;\;\;\;
\end{eqnarray}
\begin{eqnarray}
& &\;\;\;\<w_{4}',\left({\cal{Y}}_{4}(w_{1},z_{0})Y(a,z)w_{2}\right)(z_{2})w_{3}\>\nonumber\\
& &=\<w_{4}',\psi(w_{1},Y(a,z)w_{2};z_{0},z_{2})w_{3}\>\nonumber\\
& &=E_{32}
\<w_{4}',{\cal{Y}}_{1}(w_{1},z_{3}){\cal{Y}}_{2}(Y(a,z)w_{2},z_{2})w_{3}\>|_{z_{3}=z_{0}+z_{2}},
\end{eqnarray}
and
\begin{eqnarray}
& &\;\;\; \<w_{4}',\left({\cal{Y}}_{4}(Y(a,z_{4})w_{1},z_{0})w_{2}\right)(z_{2})w_{3}\>\nonumber\\
& &=E_{32}
\<w_{4}',{\cal{Y}}_{1}(Y(a,z_{4})w_{1},z_{3}){\cal{Y}}_{2}(w_{2},z_{2})w_{3}\>|_{z_{3}=
z_{0}+z_{2}}.
\end{eqnarray}
Then 
\begin{eqnarray}
& &\;\;\;\<w_{4}',\left(Y(a,z){\cal{Y}}_{4}(w_{1},z_{0})w_{2}\right)(z_{2})w_{3}\>
-\<w_{4}',\left({\cal{Y}}_{4}(w_{1},z_{0})Y(a,z)w_{2}\right)(z_{2})w_{3}\>\nonumber\\
& &=E_{32}
\<w_{4}',{\rm Res}_{z_{1}} z^{-1}\delta\left(\frac{z_{1}-z_{2}}{z}\right)
Y(a,z_{1}){\cal{Y}}_{1}(w_{1},z_{3}){\cal{Y}}_{2}(w_{2},z_{2})w_{3}\>|_{z_{3}=z_{0}+z_{2}}
\nonumber\\
& &\;\;-E_{32}
\<w_{4}',{\rm Res}_{z_{1}}z^{-1}\delta\left(\frac{z_{2}-z_{1}}{-z}\right)
{\cal{Y}}_{1}(w_{1},z_{3}){\cal{Y}}_{2}(w_{2},z_{2})Y(a,z_{1})w_{3}\>|_{z_{3}=z_{0}+z_{2}}
\nonumber\\
& &\;\;-E_{32}
\<w_{4}',{\cal{Y}}_{1}(w_{1},z_{3}){\cal{Y}}_{2}(Y(a,z)w_{2},z_{2})w_{3}\>|_{z_{3}
=z_{0}+z_{2}}\nonumber\\
& &=E_{32}{\rm Res}_{z_{1}}{\rm Res}_{z_{4}} 
z^{-1}\delta\left(\frac{z_{1}-z_{2}}{z}\right)
z_{1}^{-1}\delta\left(\frac{z_{3}+z_{4}}{z_{1}}\right)\nonumber\\
& &\;\;\;\cdot \<w_{4}',
{\cal{Y}}_{1}(Y(a,z_{4})w_{1},z_{3}){\cal{Y}}_{2}(w_{2},z_{2})w_{3}\>|_{z_{3}
=z_{0}+z_{2}}\nonumber\\
& &=E_{32}\<w_{4}',{\rm Res}_{z_{4}}z^{-1}\delta\left(\frac{z_{0}+z_{4}}{z}
\right)
{\cal{Y}}_{1}(Y(a,z_{4})w_{1},z_{3}){\cal{Y}}_{2}(w_{2},z_{2})w_{3}\>|_{z_{3}=z_{0}+z_{2}}
\nonumber\\
& &={\rm Res}_{z_{4}}z^{-1}\delta\left(\frac{z_{0}+z_{4}}{z}\right)
\<w_{4}', \left({\cal{Y}}_{4}(Y(a,z_{4})w_{1},z_{0})w_{2}\right)(z_{2})w_{3}\>.
\end{eqnarray}
This gives rise to the commutator formula which implies commutativity. Next, 
we prove associativity.
Let $k$ be a positive integer so large that 
\begin{eqnarray}\label{e3.26}
(z_{1}-z_{2})^{k}Y(a,z_{1}){\cal{Y}}_{2}(w_{2},z_{2})
=(z_{1}-z_{2})^{k}{\cal{Y}}_{2}(w_{2},z_{2})Y(a,z_{1}).
\end{eqnarray}
Since 
\begin{eqnarray}
\left({\partial \over\partial z}-{\partial \over\partial z_{2}}\right)
z^{-1}\delta\left(\frac{z_{1}-z_{2}}{z}\right)=
\left({\partial \over\partial z}-{\partial \over\partial z_{2}}\right)
z_{1}^{-1}\delta\left(\frac{z+z_{2}}{z_{1}}\right)=0,
\end{eqnarray}
we obtain
\begin{eqnarray}\label{e3.27}
& &z^{-1}\delta\left(\frac{z_{1}-z_{2}}{z}\right)
=e^{-z_{0}{\partial \over\partial z}+z_{0}{\partial \over\partial z_{2}}}
z^{-1}\delta\left(\frac{z_{1}-z_{2}}{z}\right)
=(z-z_{0})^{-1}\delta\left(\frac{z_{1}-z_{2}-z_{0}}{z-z_{0}}\right).\;\;\;\;\;
\end{eqnarray}
Similarly we have
\begin{eqnarray}\label{e3.28}
z^{-1}\delta\left(\frac{-z_{2}+z_{1}}{z}\right)
=(z-z_{0})^{-1}\delta\left(\frac{-z_{2}-z_{0}+z_{1}}{z-z_{0}}\right).
\end{eqnarray}
Then using (\ref{e3.21}), (\ref{e3.26}), (\ref{e3.27}) and (\ref{e3.28}) we obtain
\begin{eqnarray}
& &\;\;\;\;\<w_{4}',z^{k}\left(Y(a,z){\cal{Y}}_{4}(w_{1},z_{0})w_{2}\right)(z_{2})w_{3}\>
\nonumber\\
& &=E_{32}
\<w_{4}',{\rm Res}_{z_{1}} z^{-1}\delta\left(\frac{z_{1}-z_{2}}{z}\right)z^{k}
Y(a,z_{1}){\cal{Y}}_{1}(w_{1},z_{3}){\cal{Y}}_{2}(w_{2},z_{2})w_{3}\>|_{z_{3}=z_{0}+z_{2}}
\nonumber\\
& &\;\;-E_{32}
\<w_{4}',{\rm Res}_{z_{1}}z^{-1}\delta\left(\frac{z_{2}-z_{1}}{-z}\right)z^{k}
{\cal{Y}}_{1}(w_{1},z_{3}){\cal{Y}}_{2}(w_{2},z_{2})Y(a,z_{1})w_{3}\>|_{z_{3}=z_{0}+z_{2}}
\nonumber\\
& &=E_{32}
\<w_{4}',{\rm Res}_{z_{1}} z^{-1}\delta\left(\frac{z_{1}-z_{2}}{z}\right)(z_{1}-z_{2})^{k}
Y(a,z_{1}){\cal{Y}}_{1}(w_{1},z_{3}){\cal{Y}}_{2}(w_{2},z_{2})w_{3}\>_{z_{3}=z_{0}+z_{2}}\nonumber\\
& &\;\;-E_{32}
\<w_{4}',{\rm Res}_{z_{1}}z^{-1}\delta\left(\frac{z_{2}-z_{1}}{-z}\right)(z_{1}-z_{2})^{k}
{\cal{Y}}_{1}(w_{1},z_{3}){\cal{Y}}_{2}(w_{2},z_{2})Y(a,z_{1})w_{3}\>_{z_{3}=z_{0}+z_{2}}
\nonumber\\
& &={\rm Res}_{z_{1}} z^{-1}\delta\left(\frac{z_{1}-z_{2}}{z}\right)
E_{32}
\<w_{4}',(z_{1}-z_{2})^{k}
Y(a,z_{1}){\cal{Y}}_{1}(w_{1},z_{3}){\cal{Y}}_{2}(w_{2},z_{2})w_{3}\>|_{z_{3}=z_{0}+z_{2}}
\nonumber\\
& &\;\;-{\rm Res}_{z_{1}}z^{-1}\delta\left(\frac{z_{2}-z_{1}}{-z}\right)
E_{32}\<w_{4}',(z_{1}-z_{2})^{k}
{\cal{Y}}_{1}(w_{1},z_{3})Y(a,z_{1}){\cal{Y}}_{2}(w_{2},z_{2})w_{3}\>|_{z_{3}+z_{0}+z_{2}}
\nonumber\\
& &={\rm Res}_{z_{1}} (z-z_{0})^{-1}\delta\left(\frac{z_{1}-z_{2}-z_{0}}{z-z_{0}}\right)
\nonumber\\
& &\;\;\;\cdot E_{32}
\<w_{4}',(z_{1}-z_{2})^{k}
Y(a,z_{1}){\cal{Y}}_{1}(w_{1},z_{3}){\cal{Y}}_{2}(w_{2},z_{2})w_{3}\>|_{z_{3}=z_{0}+z_{2}}
\nonumber\\
& &\;\;-{\rm Res}_{z_{1}}(z-z_{0})^{-1}\delta\left(\frac{-z_{2}-z_{0}+z_{1}}{z-z_{0}}\right)
\nonumber\\
& &\;\;\;\cdot E_{32}\<w_{4}',(z_{1}-z_{2})^{k}
{\cal{Y}}_{1}(w_{1},z_{3})Y(a,z_{1}){\cal{Y}}_{2}(w_{2},z_{2})w_{3}\>|_{z_{3}=z_{0}+z_{2}}
\nonumber\\
& &={\rm Res}_{z_{1}} 
E_{32}\nonumber\\
& &\<w_{4}',(z-z_{0})^{-1}\delta\left(\frac{z_{1}-z_{3}}{z-z_{0}}\right)(z_{1}-z_{2})^{k}
Y(a,z_{1}){\cal{Y}}_{1}(w_{1},z_{3}){\cal{Y}}_{2}(w_{2},z_{2})w_{3}\>|_{z_{3}=z_{0}+z_{2}}
\nonumber\\
& &\;\;-{\rm Res}_{z_{1}}
E_{32}\nonumber\\
& &\<w_{4}',(z-z_{0})^{-1}\delta\left(\frac{-z_{3}+z_{1}}{z-z_{0}}\right)
(z_{1}-z_{2})^{k}
{\cal{Y}}_{1}(w_{1},z_{3})Y(a,z_{1}){\cal{Y}}_{2}(w_{2},z_{2})w_{3}\>|_{z_{3}=z_{0}+z_{2}}
\nonumber\\
& &={\rm Res}_{z_{1}} E_{32}
\<w_{4}',z_{1}^{-1}\delta\left(\frac{z_{2}+z}{z_{1}}\right)(z_{1}-z_{2})^{k}
{\cal{Y}}_{1}(Y(a,z-z_{0})w_{1},z_{0}+z_{2}){\cal{Y}}_{2}(w_{2},z_{2})w_{3}\>\nonumber\\
& &=\<w_{4}', z^{k}\left({\cal{Y}}_{4}(Y(a,z-z_{0})w_{1},z_{0})w_{2}\right)(z_{2})w_{3}\>.
\end{eqnarray}
This proves associativity. Combining commutativity and associativity,
we obtain the Jacobi identity.  $\;\;\;\;\Box$

Summarizing what we have proved  we obtain

\bt{t3.18} Let ${\cal{Y}}_{1}, {\cal{Y}}_{2}$ be any intertwining operators of types
$\!\left(\!\!\begin{array}{c}W_{4}\\W_{1} W_{5}\end{array}\!\!\right)$ and 
$\!\left(\!\!\begin{array}{c}W_{5}\\W_{2} W_{3}\end{array}\!\!\right)$ respectively. Then
there is a weak $V$-module $W$ and 
two intertwining operators ${\cal{Y}}_{3}, {\cal{Y}}_{4}$ of types
$\!\left(\!\!\begin{array}{c}W\\W_{1} W_{2}\end{array}\!\!\right)$ and 
$\!\left(\!\!\begin{array}{c}W_{4}\\W W_{3}\end{array}\!\!\right)$, respectively, 
such that
\begin{eqnarray}
& &\;\;\;\;E_{12}\<w_{4}',
{\cal{Y}}_{1}(w_{1},z_{1}){\cal{Y}}_{2}(w_{2},z_{2})w_{3}\>
|_{z_{1}=z_{0}+z_{2}}\nonumber\\
& &=\<w_{4}',{\cal{Y}}_{4}({\cal{Y}}_{3}(w_{1},z_{0})w_{2},z_{2})w_{3}\>.
\end{eqnarray}
for any $w_{i}\in W_{i}\; (i=1,2,3), w_{4}'\in W_{4}'$.
\et

Let $(W_{12},F_{12})$ and $ (W_{12,3},F_{12,3})$ be tensor products of
ordered pairs $(W_{1},W_{2})$ and $(W_{12},W_{3})$, respectively. Let
$(W_{23},F_{23})$ and $(W_{1,23},F_{1,23})$ be tensor products of
ordered pairs $(W_{2},W_{3})$ and $(W_{1},W_{23})$,
respectively. Taking $W_{4}=W_{1,23}$ and $W_{5}=W_{23}$, then
combining Theorem \ref{t3.18} with Proposition \ref{pu2}, we obtain
the following proposition:

\bp{p3.19} Suppose that $W$ is a weak $V$-module in the category ${\cal D}$.
 Then
there is a $V$-homomorphism $f$ from $W_{12,3}$ to $W_{1,23}$ such that
\begin{eqnarray}
& &\;\;\;\<w_{4}',fF_{12,3}(F_{12}(w_{1},z_{0})w_{2},z_{2})w_{3}\>\nonumber\\
& &=\left(E_{12}\<w_{4}',F_{1,23}(w_{1},z_{1})F_{23}(w_{2},z_{2})w_{3}\>\right)
|_{z_{1}=z_{0}+z_{2}}
\end{eqnarray}
for $w_{i}\in W_{i}\;(i=1,2,3), w_{4}'\in W_{1,23}'$.
\ep

Let $(W_{12},F_{12}(\cdot,z))$ and $(W_{12,3},F_{12,3}(\cdot,z))$ be tensor products of
ordered pairs $(W_{1},W_{2})$ and $(W_{12},W_{3})$, respectively. For any 
$w_{2}\in W_{2}, w_{3}\in W_{3}\in W_{3}$, we define 
$\phi(w_{2},w_{3};z_{1},z_{2})\in \left({\rm Hom}_{{\C}}(W_{1},W_{12,3})\right)\{z_{1},z_{2}\}$ 
as follows:
\begin{eqnarray}
& &\;\;\;\;\<w_{4}', \phi(w_{2},w_{3};z_{1},z_{2})w_{1}\>\nonumber\\
& &=E^{12}\left(\<w_{4}', F_{12,3}(F_{12}(w_{1},z_{1}-z_{2})w_{2},z_{2})w_{3}\>
\right)
\nonumber\\
& &=\sum_{i=1}^{k}z_{1}^{r_{i}}z_{2}^{s_{i}}\tilde{f}_{i}
\left({z_{2}\over z_{1}}\right),
\end{eqnarray}
for any $w_{1}\in W_{1},w_{4}'\in W_{12,3}'$. Set
\begin{eqnarray}
\phi(w_{2},w_{3};z_{1},z_{2})=\sum_{n\in {\C}}\phi_{n}(w_{2},w_{3},z_{1})z_{2}^{-n-1}
\end{eqnarray}
for any $w_{2}\in W_{2}, w_{3}\in W_{3}$. Then it follows from the
assumptions that there exists a real number $h$ such that
$\phi_{n}(w_{2},w_{3},z_{1})=0$ whenever Re$\,n>h$.

\bp{p3.20}
For $w_{2}\in W_{2},w_{3}\in W_{3},n\in {\C}$, 
$\phi_{n}(w_{2},w_{3},z_{1})\in G_{r}^{o}(W_{1},W_{12,3})$.
\ep

{\bf Proof.} For any $n\in {\C}$, we can prove as in Proposition 
\ref{p3.11} that
$\phi_{n}(w_{2},w_{3},z_{1})$ satisfies the axiom (Gr1).
Let $w_{1}\in W_{1}, w_{4}'\in W_{12,3}'$. Then
\begin{eqnarray}
& &\;\;\;\;\< w_{4}',\phi(w_{2},w_{3};z_{1},z_{2})L(-1)w_{1}\>\nonumber\\
& &=E^{12}\< w_{4}', F_{12,3}(F_{12}(L(-1)w_{1},z_{1}-z_{2})w_{2},z_{2})w_{3}\>
\nonumber\\
& &=E^{12}\left({\partial\over \partial z_{0}}
\< w_{4}', F_{12,3}(F_{12}(w_{1},z_{0})w_{2},z_{2})w_{3}\>
\right)|_{z_{0}=z_{1}-z_{2}}\nonumber\\
& &={\partial\over \partial z_{1}}E^{12}
\< w_{4}', F_{12,3}(F_{12}(w_{1},z_{0})w_{2},z_{2})w_{3}\>
|_{z_{0}=z_{1}-z_{2}}\nonumber\\
& &={\partial\over \partial z_{1}}\< w_{4}',\phi(w_{2},w_{3};z_{1},z_{2})w_{1}\>.
\end{eqnarray}
This implies that $\phi_{n}(w_{2},w_{3},z_{1})$ satisfies the axiom (Gr2).
For any $a\in V, n\in {\C}$, let $k$ be a positive integer such that
\begin{eqnarray}
& &\phi_{h}(w_{2},w_{3},z_{1})=0\;\;\;\mbox{ whenever Re}\,h>k+n;\\
& &z^{k}Y(a,z)F_{12}(w_{1},z_{2})w_{3}=z^{k}F_{12}(Y(a,z-z_{2})w_{1},z_{2})w_{3};\\
& &z^{k}Y(a,z)F_{12,3}(w,z_{2})w_{3}=z^{k}F_{12,3}(Y(a,z-z_{2})w,z_{2})w_{3},
\end{eqnarray}
for any $w_{1}\in W_{1}, w\in W_{12}$.  Then we obtain
\begin{eqnarray}
& &\;\;\;z^{3k}\< w_{4}',Y(a,z)\phi_{n}(w_{2},w_{3};z_{1})w_{1}\>\nonumber\\
& &={\rm Res}_{z_{2}}z_{2}^{n}
z^{3k}\< w_{4}',Y(a,z)\phi(w_{2},w_{3};z_{1},z_{2})w_{1}\>\nonumber\\
& &={\rm Res}_{z_{2}}\sum_{i=0}^{2k}{2k\choose i}(z-z_{2})^{2k-i}z_{2}^{n+i}z^{k}
\< w_{4}',Y(a,z)\phi(w_{2},w_{3};z_{1},z_{2})w_{1}\>\nonumber\\
& &={\rm Res}_{z_{2}}\sum_{i=0}^{k}{2k\choose i}(z-z_{2})^{2k-i}z_{2}^{n+i}z^{k}
\< w_{4}',Y(a,z)\phi(w_{2},w_{3};z_{1},z_{2})w_{1}\>\nonumber\\
& &={\rm Res}_{z_{2}}\sum_{i=0}^{k}{2k\choose i}(z-z_{2})^{2k-i}z_{2}^{n+i}z^{k}\cdot \nonumber\\
& &\;\;\;\cdot E^{12}\< w_{4}', F_{12,3}(Y(a,z-z_{2})
F_{12}(w_{1},z_{0})w_{2},z_{2})w_{3}\>|_{z_{0}=z_{1}-z_{2}}\nonumber\\
& &={\rm Res}_{z_{2}}\sum_{i=0}^{k}{2k\choose i}(z-z_{2})^{2k-i}z_{2}^{n+i}z^{k}\cdot \nonumber\\
& &\;\;\;\cdot E^{12}\< w_{4}', 
F_{12,3}(F_{12}(Y(a,z-z_{2}-z_{0})w_{1},z_{0})w_{2},z_{2})
w_{3}\>|_{z_{0}=z_{1}-z_{2}}\nonumber\\
& &={\rm Res}_{z_{2}}\sum_{i=0}^{k}{2k\choose i}(z-z_{2})^{2k-i}z_{2}^{n+i}z^{k}\cdot \nonumber\\
& &\;\;\;\cdot E^{12}\< w_{4}', 
F_{12,3}(F_{12}(Y(a,z-z_{1})w_{1},z_{0})w_{2},z_{2})
w_{3}\>|_{z_{0}=z_{1}-z_{2}}\nonumber\\
& &={\rm Res}_{z_{2}}\sum_{i=0}^{2k}{2k\choose i}(z-z_{2})^{2k-i}z_{2}^{n+i}z^{k}\cdot \nonumber\\
& &\;\;\;\cdot E^{12}\< w_{4}', 
F_{12,3}(F_{12}(Y(a,z-z_{1})w_{1},z_{0})w_{2},z_{2})
w_{3}\>|_{z_{0}=z_{1}-z_{2}}\nonumber\\
& &={\rm Res}_{z_{2}}z^{3k}z_{2}^{n}
E^{12}\< w_{4}',F_{12,3}(F_{12}(Y(a,z-z_{1})w_{1},z_{0})w_{2},z_{2})w_{3}\>
|_{z_{0}=z_{1}-z_{2}}\nonumber\\
& &={\rm Res}_{z_{2}}z^{3k}z_{2}^{n}\< w_{4}',\phi(w_{2},w_{3};z_{1},z_{2})Y(a,z-z_{1})w_{1}\>
\nonumber\\
& &=z^{3k}\< w_{4}',\phi(w_{2},w_{3};z_{1},z_{2})Y(a,z-z_{1})w_{1}\>.
\end{eqnarray}
Then $\phi_{n}(w_{2},w_{3},z_{1})$ satisfies the axiom (Gr3). Thus 
$\phi_{n}(w_{2},w_{3},z_{1})\in G_{r}(W_{1},W_{12,3})$. $\;\;\;\;\Box$

As in Lemma \ref{l3.2} we have

\begin{lem}\label{l3.21}
For $a\in V, m\in {\Z}, w_{2}\in W_{2},w_{3}\in W_{3}$ we have
\begin{eqnarray}
a_{m}\circ_{r} \phi(w_{2},w_{3};z_{1},z_{2})
=\phi(w_{2},a_{m}w_{3};z_{1},z_{2})+\sum_{j=0}^{\infty}{m\choose j}
z_{2}^{m-j}\phi(a_{j}w_{2},w_{3};z_{1},z_{2}).
\end{eqnarray}
\end{lem}

{\bf Proof.} For any $w_{1}\in W_{1},w_{4}'\in W_{12,3}'$, we have
\begin{eqnarray}
& &\;\;\;\;\< w_{4}',\left(a_{m}\circ_{r} \phi(w_{2},w_{3};z_{1},z_{2})\right)w_{1}\>\nonumber\\
& &=\< w_{4}',a_{m}\phi(w_{2},w_{3};z_{1},z_{2})w_{1}\>
-\sum_{i=0}^{\infty}{m\choose i}\< w_{4}',z_{1}^{m-i}\phi(w_{2},w_{3};z_{1},z_{2})a_{i}w_{1}\>
\nonumber\\
& &=E^{12}\< w_{4}',a_{m}F_{12,3}(F_{12}(w_{1},z_{0})w_{2},z_{2})w_{3}\>
|_{z_{0}=z_{1}-z_{2}}\nonumber\\
& &\;\;\;-\sum_{i=0}^{\infty}{m\choose i}E^{12}\< w_{4}',z_{1}^{m-i}
F_{12,3}(F_{12}(a_{i}w_{1},z_{0})w_{2},z_{2})w_{3}\>
|_{z_{0}=z_{1}-z_{2}}\nonumber\\
& &=E^{12}\< w_{4}',F_{12,3}(F_{12}(w_{1},z_{0})w_{2},z_{2})a_{m}w_{3}\>
|_{z_{0}=z_{1}-z_{2}}\nonumber\\
& &\;\;\;+\sum_{j=0}^{\infty}{m\choose j}E^{12}\< w_{4}',z_{2}^{m-j}
F_{12,3}(a_{j}F_{12}(w_{1},z_{0})w_{2},z_{2})w_{3}\>
|_{z_{0}=z_{1}-z_{2}}\nonumber\\
& &\;\;\;-\sum_{i=0}^{\infty}{m\choose i}E^{12}\< w_{4}',z_{1}^{m-i}
F_{12,3}(F_{12}(a_{i}w_{1},z_{0})w_{2},z_{2})w_{3}\>
|_{z_{0}=z_{1}-z_{2}}\nonumber\\
& &=\< w_{4}',\phi(w_{2},a_{m}w_{3};z_{1},z_{2})w_{1}\>\nonumber\\
& &\;\;\;+\sum_{j=0}^{\infty}{m\choose j}E^{12}\< w_{4}',z_{2}^{m-j}
F_{12,3}(F_{12}(w_{1},z_{0})a_{j}w_{2},z_{2})w_{3}\>
|_{z_{0}=z_{1}-z_{2}}\nonumber\\
& &\;\;\;+\sum_{j=0}^{\infty}\sum_{i=0}^{\infty}{m\choose j}{j\choose i}
E^{12}\< w_{4}',z_{2}^{m-j}z_{0}^{j-i}
F_{12,3}(F_{12}(a_{i}w_{1},z_{0})w_{2},z_{2})w_{3}\>
|_{z_{0}=z_{1}-z_{2}}\nonumber\\
& &\;\;\;-\sum_{i=0}^{\infty}{m\choose i}E^{12}\< w_{4}',z_{1}^{m-i}
F_{12,3}(F_{12}(a_{i}w_{1},z_{0})w_{2},z_{2})w_{3}\>
|_{z_{0}=z_{1}-z_{2}}\nonumber\\
& &=\< w_{4}',\phi(w_{2},a_{m}w_{3};z_{1},z_{2})w_{1}\>
+\sum_{j=0}^{\infty}{m\choose j}\< w_{4}',z_{2}^{m-j}\phi(a_{j}w_{2},w_{3};z_{1},z_{2})w_{1}\>.
\end{eqnarray}
Then the lemma follows immediately. $\;\;\;\;\Box$

It follows from Lemma \ref{l3.21} that the space $M$ linearly spanned by all 
$\phi_{n}(w_{2},w_{3},z_{1})$ is a weak submodule of $G_{r}^{o}(W_{1},W_{12,3})$.
Then we obtain a natural intertwining operator $I_{1}(\cdot,z)$ of type
$\left(\!\begin{array}{c}W_{12,3}\\W_{1} M\end{array}\!\right)$ such that
\begin{eqnarray}
I_{1}(w_{1},z)A=A(z)w_{1}\;\;\;\mbox{ for any } w_{1}\in W_{1}, A\in M.
\end{eqnarray}
Next we define a linear map $I_{2}(\cdot,z)$ from $W_{2}$ to 
$\left({\rm Hom}_{{\C}}(W_{3},M)\right)\{z\}$ as follows:
\begin{eqnarray}
\left(I_{2}(w_{2},z_{2})w_{3}\right)(z_{1})=\phi(w_{2},w_{3};z_{1},z_{2})
\end{eqnarray}
for any $w_{2}\in W_{2},w_{3}\in W_{3}$.

\bp{p3.22}
The  map $I_{2}(\cdot,z)$ is an intertwining operator of type
$\left(\!\!\begin{array}{c}M\\W_{2} W_{3}\end{array}\!\!\right)$.
\ep

{\bf Proof.} Since there exists a real number $h$ such that $\phi_{n}(w_{2},w_{3},z_{1})=0$ whenever
$n>h$, $I_{2}(\cdot,z)$ satisfies the axiom (I1). Let 
$w_{i}\in W_{i}\;(i=1,2,3), w_{4}'\in W_{12,3}'$. Then we have
\begin{eqnarray}
& &\;\;\;\< w_{4}', \phi(L(-1)w_{2},w_{3};z_{1},z_{2})w_{1}\>\nonumber\\
& &=E^{12}\< w_{4}', F_{12,3}(F_{12}(w_{1},z_{0})L(-1)w_{2},z_{2})w_{3}\>
|_{z_{0}=z_{1}-z_{2}}\nonumber\\
& &=E^{12}\< w_{4}', F_{12,3}(L(-1)F_{12}(w_{1},z_{0})w_{2},z_{2})w_{3}\>
|_{z_{0}=z_{1}-z_{2}}\nonumber\\
& &\;\;\;-E^{12}\left(\< w_{4}', {\partial\over\partial z_{0}}
F_{12,3}(F_{12}(w_{1},z_{0})w_{2},z_{2})w_{3}\>
\right)|_{z_{0}=z_{1}-z_{2}}\nonumber\\
& &=E^{12}\left(\< w_{4}', {\partial\over\partial z_{2}}
F_{12,3}(F_{12}(w_{1},z_{0})w_{2},z_{2})w_{3}\>
\right)|_{z_{0}=z_{1}-z_{2}}\nonumber\\
& &\;\;\;-E^{12}\left(\< w_{4}', {\partial\over\partial z_{0}}
F_{12,3}(F_{12}(w_{1},z_{0})w_{2},z_{2})w_{3}\>
\right)|_{z_{0}=z_{1}-z_{2}}\nonumber\\
& &={\partial\over\partial z_{2}}E^{12}\< w_{4}', 
F_{12,3}(F_{12}(w_{1},z_{0})w_{2},z_{2})w_{3}\>
|_{z_{0}=z_{1}-z_{2}}\nonumber\\
& &={\partial\over\partial z_{2}}
\< w_{4}', \phi(w_{2},w_{3};z_{1},z_{2})w_{1}\>.
\end{eqnarray}
This gives the $L(-1)$-derivative property. To prove the Jacobi
identity, as in Proposition \ref{p3.3x} we prove the commutator
formula and the associativity.  For $a\in V$, we have
\begin{eqnarray}
& &\;\;\;\<w_{4}',(Y(a,z)I_{2}(w_{2},z_{2})w_{3})(z_{1})w_{1}\>
-\<w_{4}',(I_{2}(w_{2},z_{2})Y(a,z)w_{3})(z_{1})w_{1}\>\nonumber\\
& &=\<w_{4}',Y(a,z)(I_{2}(w_{2},z_{2})w_{3})(z_{1})w_{1}\>\nonumber\\
& &\;\;\;-{\rm Res}_{z_{3}}z^{-1}\delta\left(\frac{z_{1}+z_{3}}{z}\right)
\<w_{4}',(I_{2}(w_{2},z_{2})w_{3})(z_{1})Y(a,z_{3})w_{1}\>\nonumber\\
& &\;\;\;-\<w_{4}',\phi(w_{2},Y(a,z)w_{3};z_{1},z_{2})w_{1}\>\nonumber\\
& &=\<w_{4}',Y(a,z)\phi(w_{2},w_{3};z_{1},z_{2})w_{1}\>\nonumber\\
& &\;\;\;-{\rm Res}_{z_{3}}z^{-1}\delta\left(\frac{z_{1}+z_{3}}{z}\right)
\<w_{4}',\phi(w_{2},w_{3};z_{1},z_{2})Y(a,z_{3})w_{1}\>\nonumber\\
& &\;\;\;-E^{12}(\<w_{4}',F_{12,3}(F_{12}(w_{1},z_{0})w_{2},z_{2})
Y(a,z)w_{3}\>|_{z_{0}=z_{1}-z_{2}}.
\nonumber\\
& &=E^{12}\<w_{4}',Y(a,z)F_{12,3}(F_{12}(w_{1},z_{0})w_{2},z_{2})w_{3}\>
|_{z_{0}=z_{1}-z_{2}}\nonumber\\
& &\;\;\;-{\rm Res}_{z_{3}}z^{-1}\delta\left(\frac{z_{1}+z_{3}}{z}\right)
E^{12}\<w_{4}',F_{12,3}(F_{12}(Y(a,z_{3})w_{1},z_{0})w_{2},z_{2})w_{3}\>
|_{z_{0}=z_{1}-z_{2}}\nonumber\\
& &\;\;\;-E^{12}\<w_{4}',F_{12,3}(F_{12}(w_{1},z_{0})w_{2},z_{2})
Y(a,z)w_{3}\>|_{z_{0}=z_{1}-z_{2}}
\nonumber\\
& &=E^{12}\<w_{4}',{\rm Res}_{z_{3}}z^{-1}\delta\left(\frac{z_{2}+z_{3}}{z}\right)
F_{12,3}(Y(a,z_{3})F_{12}(w_{1},z_{0})w_{2},z_{2})w_{3}\>|_{z_{0}=z_{1}-z_{2}}
\nonumber\\
& &\;\;\;-{\rm Res}_{z_{3}}z^{-1}\delta\left(\frac{z_{1}+z_{3}}{z}\right)
E^{12}\<w_{4}',F_{12,3}(F_{12}(Y(a,z_{3})w_{1},z_{0})w_{2},z_{2})w_{3}\>
|_{z_{0}=z_{1}-z_{2}}\nonumber\\
& &=E^{12}\<w_{4}',{\rm Res}_{z_{3}}z^{-1}\delta\left(\frac{z_{2}+z_{3}}{z}\right)
F_{12,3}(F_{12}(w_{1},z_{0})Y(a,z_{3})w_{2},z_{2})w_{3}\>|_{z_{0}=z_{1}-z_{2}}
\nonumber\\
& &\;\;\;+{\rm Res}_{z_{3}}z^{-1}\delta\left(\frac{z_{2}+z_{3}}{z}\right)\cdot\nonumber\\
& &
\cdot E^{12}\<w_{4}',{\rm Res}_{z_{4}}
z_{0}^{-1}\delta\left(\frac{z_{3}-z_{4}}{z_{0}}\right)
F_{12,3}(F_{12}(Y(a,z_{4})w_{1},z_{0})w_{2},z_{2})w_{3}\>|_{z_{0}=z_{1}-z_{2}}
\nonumber\\
& &\;\;\;-{\rm Res}_{z_{3}}z^{-1}\delta\left(\frac{z_{1}+z_{3}}{z}\right)
E^{12}\<w_{4}',F_{12,3}(F_{12}(Y(a,z_{3})w_{1},z_{0})w_{2},z_{2})w_{3}\>
|_{z_{0}=z_{1}-z_{2}}\nonumber\\
& &={\rm Res}_{z_{3}}z^{-1}\delta\left(\frac{z_{2}+z_{3}}{z}\right)
\<w_{4}',\phi(Y(a,z_{3})w_{2},w_{3};z_{1},z_{2})w_{1}\>\nonumber\\
& &\;\;\;+E^{12}\<w_{4}',{\rm Res}_{z_{4}}
z^{-1}\delta\left(\frac{z_{2}+z_{0}+z_{4}}{z}\right)
F_{12,3}(F_{12}(Y(a,z_{4})w_{1},z_{0})w_{2},z_{2})w_{3}\>|_{z_{0}=z_{1}-z_{2}}
\nonumber\\
& &\;\;\;-{\rm Res}_{z_{3}}z^{-1}\delta\left(\frac{z_{1}+z_{3}}{z}\right)
E^{12}\<w_{4}',F_{12,3}(F_{12}(Y(a,z_{3})w_{1},z_{0})w_{2},z_{2})w_{3}\>
|_{z_{0}=z_{1}-z_{2}}\nonumber\\
& &={\rm Res}_{z_{3}}z^{-1}\delta\left(\frac{z_{2}+z_{3}}{z}\right)
\<w_{4}',\phi(Y(a,z_{3})w_{2},w_{3};z_{1},z_{2})w_{1}\>\nonumber\\
& &\;\;\;+E^{12}\<w_{4}',{\rm Res}_{z_{4}}
z^{-1}\delta\left(\frac{z_{1}+z_{4}}{z}\right)
F_{12,3}(F_{12}(Y(a,z_{4})w_{1},z_{0})w_{2},z_{2})w_{3}\>|_{z_{0}=z_{1}-z_{2}}
\nonumber\\
& &\;\;\;-{\rm Res}_{z_{3}}z^{-1}\delta\left(\frac{z_{1}+z_{3}}{z}\right)
E^{12}\<w_{4}',F_{12,3}(F_{12}(Y(a,z_{3})w_{1},z_{0})w_{2},z_{2})w_{3}\>
|_{z_{0}=z_{1}-z_{2}}\nonumber\\
& &={\rm Res}_{z_{3}}z^{-1}\delta\left(\frac{z_{2}+z_{3}}{z}\right)
\<w_{4}',\phi(Y(a,z_{3})w_{2},w_{3};z_{1},z_{2})w_{1}\>\nonumber\\
& &={\rm Res}_{z_{3}}z^{-1}\delta\left(\frac{z_{2}+z_{3}}{z}\right)
\<w_{4}',(I_{2}(Y(a,z_{3})w_{2},z_{2})w_{3})(z_{1})w_{1}\>.
\end{eqnarray}
This proves the commutator formula. Let $k$ be a positive integer such that
\begin{eqnarray}
& &z_{1}^{k}Y(a,z_{1})F_{12,3}(u,z_{2})w_{3}=z_{1}^{k}F_{12,3}(Y(a,z_{1}-z_{2})u,z_{2})w_{3};\\
& &z_{1}^{k}Y(a,z_{1})F_{12}(v,z_{2})w_{3}=z_{1}^{k}F_{12}(Y(a,z_{1}-z_{2})v,z_{2})w_{3}
\end{eqnarray}
for any $u\in W_{12},v\in W_{1}$. Then we obtain
\begin{eqnarray}
& &\;\;\;z^{k}\<w_{4}',(Y(a,z)I_{2}(w_{2},z_{2})w_{3})(z_{1})w_{1}\>\nonumber\\
& &=E^{12}\<w_{4}',z^{k}Y(a,z)F_{12,3}(F_{12}(w_{1},z_{0})w_{2},z_{2})w_{3}\>
|_{z_{0}=z_{1}-z_{2}}\nonumber\\
& &\;\;\;-{\rm Res}_{z_{3}}z^{-1}\delta\left(\frac{z_{1}+z_{3}}{z}\right)z^{k}
E^{12}\<w_{4}',F_{12,3}(F_{12}(Y(a,z_{3})w_{1},z_{0})w_{2},z_{2})w_{3}\>
|_{z_{0}=z_{1}-z_{2}}\nonumber\\
& &=E^{12}\<w_{4}',
z^{k}F_{12,3}(Y(a,z-z_{2})F_{12}(w_{1},z_{0})w_{2},z_{2})w_{3}\>|_{z_{0}=z_{1}-z_{2}}\nonumber\\
& &\;\;\;-{\rm Res}_{z_{3}}z^{-1}\delta\left(\frac{z_{1}+z_{3}}{z}\right)z^{k}
E^{12}\<w_{4}',F_{12,3}(F_{12}(Y(a,z_{3})w_{1},z_{0})w_{2},z_{2})w_{3}\>
|_{z_{0}=z_{1}-z_{2}}\nonumber\\
& &=E^{12}\<w_{4}',
z^{k}F_{12,3}(F_{12}(w_{1},z_{0})Y(a,z-z_{2})w_{2},z_{2})w_{3}\>|_{z_{0}=z_{1}-z_{2}}\nonumber\\
& &\;\;\;+{\rm Res}_{z_{4}}E^{12}z_{0}^{-1}\delta\left(\frac{z-z_{2}-z_{4}}{z_{0}}\right)
\<w_{4}',
z^{k}F_{12,3}(F_{12}(Y(a,z_{4})w_{1},z_{0})w_{2},z_{2})w_{3}\>|_{z_{0}=z_{1}-z_{2}}\nonumber\\
& &\;\;\;-{\rm Res}_{z_{3}}z^{-1}\delta\left(\frac{z_{1}+z_{3}}{z}\right)z^{k}
E^{12}\<w_{4}',F_{12,3}(F_{12}(Y(a,z_{3})w_{1},z_{0})w_{2},z_{2})w_{3}\>
|_{z_{0}=z_{1}-z_{2}}\nonumber\\
& &=E^{12}\<w_{4}',
z^{k}F_{12,3}(F_{12}(w_{1},z_{0})Y(a,z-z_{2})w_{2},z_{2})w_{3}\>|_{z_{0}=z_{1}-z_{2}}\nonumber\\
& &\;\;\;+{\rm Res}_{z_{4}}E^{12}z^{-1}\delta\left(\frac{z_{0}+z_{2}+z_{4}}{z}\right)
\<w_{4}',
z^{k}F_{12,3}(F_{12}(Y(a,z_{4})w_{1},z_{0})w_{2},z_{2})w_{3}\>|_{z_{0}=z_{1}-z_{2}}\nonumber\\
& &\;\;\;-{\rm Res}_{z_{3}}z^{-1}\delta\left(\frac{z_{1}+z_{3}}{z}\right)z^{k}
E^{12}\<w_{4}',F_{12,3}(F_{12}(Y(a,z_{3})w_{1},z_{0})w_{2},z_{2})w_{3}\>
|_{z_{0}=z_{1}-z_{2}}\nonumber\\
& &=E^{12}\<w_{4}',
z^{k}F_{12,3}(F_{12}(w_{1},z_{0})Y(a,z-z_{2})w_{2},z_{2})w_{3}\>|_{z_{0}=z_{1}-z_{2}}\nonumber\\
& &=z^{k}\<w_{4}',(I_{2}(Y(a,z-z_{2})w_{2},z_{2})w_{3})(z_{1})w_{1}\>.
\end{eqnarray}
This proves the associativity. Therefore $I_{2}(\cdot,z)$ is an intertwining operator.$\;\;\;\;\Box$

By Propositions \ref{pu1} and \ref{p3.22} we obtain

\bp{p3.23} Assume that $M$ is in category ${\cal D}.$ Then 
there is a $V$-homomorphism $g$ from $W_{1,23}$ to $W_{12,3}$ such that
\begin{eqnarray}
& &\;\;\;E^{12}\<w_{4}',F_{12,3}(F_{12}(w_{1},z_{1}-z_{2})w_{2},z_{2})w_{3}\>\nonumber\\
& &=\<w_{4}',gF_{1,23}(w_{1},z_{1})F_{23}(w_{2},z_{2})w_{3}\>
\end{eqnarray}
for $w_{i}\in W_{i}\;(i=1,2,3), w_{4}'\in W_{12,3}'$.
\ep

\bt{t3.24} Let $V$ be a vertex operator algebra and ${\cal D}$ equals to
either ${\cal C}(V)$ or ${\cal C}^1(V)$ and $W_1,W_2,W_3$ weak modules in
${\cal D}.$ Assume that the tensor products 
$W_{12}=W_1\boxtimes_{\cal D}W_2,$ $W_{23}=W_2\boxtimes_{\cal D}W_3,$ 
$W_{12,3}=W_{12}\boxtimes_{\cal D}W_3,$ $W_{1,23}=W_1\boxtimes_{\cal D}W_{23}$
exist such that the convergence
and extension property for products and iterates hold. Also assume that
$W$ and $M$ are in ${\cal D}.$ Then
the $V$-homomorphisms $f$ and $g$ given by Propositions \ref{p3.19} and \ref{p3.23} are isomorphisms between $W_{1,23}$ and $W_{12,3}.$ That is, the
tensor products for $W_1,W_2,W_3$ are associative.
\et

{\bf Proof.} Let $g^*: W_{12,3}'\to W_{1,23}'$ be the adjoint map of 
$g.$ Then by Proposition \ref{p3.23} we have 
$$E^{12}\<w_{4}',F_{12,3}(F_{12}(w_{1},z_{1}-z_{2})w_{2},z_{2})w_{3}\>
=\<g^*w_{4}',F_{1,23}(w_{1},z_{1})F_{23}(w_{2},z_{2})w_{3}\>$$
for $w_i\in W_i$ and $w_4'\in W_{12,3}'.$ In Proposition \ref{p3.19}
replacing $w_4'$ by $g^*w_4'$ we obtain
$$\<g^*w_{4}',fF_{12,3}(F_{12}(w_{1},z_{1}-z_{2})w_{2},z_{2})w_{3}\>
=E_{12}\<g^*w_{4}',F_{1,23}(w_{1},z_{1})F_{23}(w_{2},z_{2})w_{3}\>.$$
Thus
$$E_{12}E^{12}\<w_{4}',F_{12,3}(F_{12}(w_{1},z_{1}-z_{2})w_{2},z_{2})w_{3}\>
=\<w_{4}',gfF_{12,3}(F_{12}(w_{1},z_{1}-z_{2})w_{2},z_{2})w_{3}\>.$$
It is obvious that the operator $E_{12}E^{12}$ is identity as
the analytic continuation of a multi-valued complex function is unique.
This implies that 
$$\<w_{4}',F_{12,3}(F_{12}(w_{1},z_{1}-z_{2})w_{2},z_{2})w_{3}\>
=\<w_{4}',gfF_{12,3}(F_{12}(w_{1},z_{1}-z_{2})w_{2},z_{2})w_{3}\>.$$
Since $W_{12,3}$ is linearly spanned by
the coefficients of $F_{12,3}(F_{12}(w_{1},z_{0})w_{2},z_{2})w_{3}$ for 
$w_{i}\in W_{i}\;(i=1,2,3)$ (see  Lemma \ref{l2.3}) and $w_i, w_4'$ are 
arbitrary  we conclude that $fg=1$. Similarly, $gf=1.$
$\;\;\;\;\Box$

Recall Remark \ref{r2.3}. 
\bt{t3.25} Let $V$ be a vertex operator algebra such that (1) ${\cal C}(V)$
is semi-simple, (2) ${\cal C}(V)$ has only finitely many simple objects,
(3) Fusion rules among any three modules are finite, (3) Any finitely generated
weak module is an ordinary module, 
(4)  The convergence
and extension property for products and iterates hold. Then ${\cal C}(V)$ is
a tensor category, that is the tensor product of any two modules exists
and the tensor product is associative.
\et

{\bf Proof.} The existence of tensor product of any two modules is given
in [HL1]-[HL4] and [LL] (cf. Remark \ref{r2.3}). By Theorem \ref{t3.25} it remains
to show that both $W$ and $M$ are in ${\cal C}(V).$ By assumption (3), 
$G_l^0(W_1,W_2)$ and $G_r^0(W_1,W_2)$ are ordinary $V$-modules. So are
$W$ and $M.$

Theorem \ref{t3.25} has been proved in [H1] using Huang-Lepowsky's approach to
the tensor product.


\begin{thebibliography}{xxxx}

\bibitem[BPZ]{BPZ}
A. Belavin, A. M. Polyakov, A. B. Zamolodchikov, Infinite conformal
symmetries in two-dimensional quantum field theory, {\it Nucl. Phys. B}
{\bf 241} (1984), 333-380.

\bibitem[B]{B}
R. E. Borcherds, Vertex algebras, Kac-Moody algebras, and the Monster,
{\it Proc. Natl. Acad. Sci. USA} {\bf 83} (1986), 3068-3071.

\bibitem[BF]{BF} D. Bernard and G. Felder, 
Fock representations and BRST cohomology in $SL(2)$ current algebra,
{\em Commun. Math. Phys.} {\bf 127} (1990), 145-168.

\bibitem[DL]{DL}
C. Dong and J. Lepowsky, Generalized Vertex Algebras and Relative
Vertex Operators, Progress in Math., {\bf Vol. 112}, Birkh\"auser,
Boston, 1993.

\bibitem[DLM1]{DLM1}
C. Dong, H. Li and G. Mason, Twisted representations of vertex operator algebras, preprint,
q-alg/9509005.

\bibitem[DLM2]{DLM2}
C. Dong, H. Li and G. Mason, Regularity of rational vertex operator algebras, preprint, 
q-alg/9508018.

\bibitem[Fe]{Fe}
G. Felder, BRST approach to minimal models, {\em Nucl. Phys.} {\bf B317} (1989),
215-236.

\bibitem[Fi]{Fi}
M. Finkelberg, Fusion categories, Ph.D. thesis, Harvard University, 1993.


\bibitem[FHL]{FHL}
I. Frenkel, Y.-Z. Huang and J. Lepowsky, On axiomatic approaches to
vertex operator algebras and modules, preprint, 1989; Memoirs Amer. Math.
Soc. {\bf 104}, 1993.

\bibitem[FLM]{FLM}
I. Frenkel, J. Lepowsky and A. Meurman, {\it Vertex Operator Algebras
and the Monster}, Pure and Appl. Math., {\bf Vol. 134}, Academic Press,
Boston, 1988.

\bibitem[KZ]{KZ}
V. G. Knizhnik and A. B. Zamolodchikov,
Current algebra and Wess-Zumino model in two dimensions, {\em Nucl. Phys.} {\bf B247} (1984),
83-103.

\bibitem[H1]{H1}
Y.-Z. Huang, A theory of tensor products for module categories for a vertex operator algebra, 
IV, {\em J. Pure Appl. Alg.} {\bf 100} (1995), 173-216.

\bibitem[H2]{H2}
Y.-Z. Huang, Virasoro vertex operator algebras, the (nonmeromorphic) operator product
expansion and the tensor product theory, preprint.

\bibitem[HL0]{HL0}
Y.-Z. Huang and J. Lepowsky, Toward a theory of tensor product for
representations 
for a vertex operator algebra, in {\it Proc. 20th
International Conference on Differential Geometric Methods in
Theoretical Physics, New York, 1991,} ed. S. Catto and A. Rocha, World
Scientific, Singapore, 1992, {\bf Vol. 1}, 344-354.

\bibitem[HL1]{HL1}
Y.-Z. Huang and J. Lepowsky, A theory of tensor product for
module category of a vertex operator algebra, I, {\em Selecta Mathematica},
to appear.

\bibitem[HL2]{HL2}
Y.-Z. Huang and J. Lepowsky, A theory of tensor product for
module category of a vertex operator algebra, II, {\em Selecta Mathematica},
to appear.

\bibitem[HL3]{HL3}
Y.-Z. Huang and J. Lepowsky, A theory of tensor product for module category
of a vertex operator algebra, III, {\em J. Pure Appl. Alg.} 
{\bf 100} (1995), 141-171.

\bibitem[HL4]{HL4}
Y.-Z. Huang and J. Lepowsky, Tensor products of 
modules for a vertex operator algebra and vertex tensor categories, 
in: {\em Lie Theory and Geometry, in honor of Bertram Kostant,} ed. R. Brylinski, 
J.-L. Brylinski, V. Guillemin, V. Kac, Birkh\"auser, 1994, 349-383.


\bibitem[KL0]{KL0}
D. Kazhdan and G. Lusztig, Affine Lie algebras and quantum groups,
{\it Duke Math. J.} {\bf 2} (1991), 21-29.

\bibitem[KL1]{KL1}
D. Kazhdan and G. Lusztig, Tensor structures arising from affine Lie
algebras, I, {\it J. Amer. Math. Soci.} {\bf Vol. 64} (1993), 905-9047

\bibitem[KL2]{KL2}
D. Kazhdan and G. Lusztig, Tensor structures arising from affine Lie
algebras, II, {\it J. Amer. Math. Soci.} {\bf Vol. 64} (1993), 9047-1011.

\bibitem[L1]{L1}
H.-S. Li, Local systems of vertex operators, vertex superalgebras and
modules, {\em J. Pure Appl. Alg.}, to appear.

\bibitem[L2]{L2}
H.-S. Li, Representation theory and tensor product theory for vertex
operator algebras, Ph.D. thesis, Rutgers University, 1994.

\bibitem[LL]{LL}
J. Lepowsky and H.-S. Li, An analogue of ``Hom''-functor and a generalized nuclear 
democracy theorem, to appear.


\bibitem[MS]{MS}
G. Moore and N. Seiberg, Classical and quantum conformal field theory,
{\it Comm. Math. Phys.} {\bf 123} (1989), 177-254.

\bibitem[TK]{TK}
A. Tsuchiya and Y. Kanie, Vertex operators in conformal field theory
on ${\bf P}^{1}$ and monodromy representations of braid group, in:
{\it Conformal Field Theory and Solvable Lattice Models, Advanced
Studies in Pure Math.}, {\bf Vol. 16}, Kinokuniya Company Ltd., Tokyo,
1988, 297-372.


\bibitem[Z]{Z}
Y.-C. Zhu, Vertex operator algebras, elliptic functions and modular
forms, Ph.D. thesis, Yale University, 1990.

\end{thebibliography}
\end{document}